\documentclass[final]{siamltex}
\usepackage{amsmath, amssymb, amsfonts}
\usepackage{graphicx,color}
\usepackage{diagbox}
\usepackage[ruled]{algorithm2e}
\usepackage{algorithmic}
\usepackage{enumitem}
\usepackage{multirow}
\usepackage{subfig}
\usepackage{comment}
\usepackage{graphicx}
\usepackage{epstopdf}
\usepackage{booktabs}
\usepackage{placeins}
\usepackage{multirow}

\newcommand{\bvec}[1]{\mathbf{#1}}

\newcommand{\mc}[1]{\mathcal{#1}}

\newcommand{\vn}{\bvec{n}}

\newcommand{\ud}{\,\mathrm{d}}

\newcommand{\abs}[1]{\lvert#1\rvert}

\newcommand{\sign}{\text{sign}}

\newcommand{\Or}{\mathcal{O}}

\newcommand{\I}{\imath} 
\newcommand{\RR}{\mathbb{R}}
\newcommand{\CC}{\mathbb{C}}

\newcommand{\jmp}[1]{\jl#1\jr}
\newcommand{\al}{\{\hspace{-3.5pt}\{}
\newcommand{\ar}{\}\hspace{-3.5pt}\}}
\newcommand{\avg}[1]{\al#1\ar}
\newcommand{\jl}{[\![}
\newcommand{\jr}{]\!]}
\newcommand{\VN}{\mathbb V_N}

\newcommand{\vpos}{\bvec{x}}

\newcommand{\LL}[1]{\textcolor{red}{[LL:#1]}}
\newcommand{\YL}[1]{\textcolor{blue}{[YL:#1]}}

\title{Globally Constructed Adaptive Local Basis Set for \\
Spectral Projectors of \\ Second Order Differential Operators}

\author{
    Yingzhou Li\thanks{Department of Mathematics, Duke University, Durham,
    NC 27708. Email: 
    \texttt{yingzhou.li@duke.edu}}
    \and
    Lin Lin\thanks{Department of Mathematics,
    University of California, Berkeley, Berkeley, CA 94720 and Computational
    Research Division, Lawrence Berkeley National Laboratory, Berkeley, CA
    94720. Email: \texttt{linlin@math.berkeley.edu}}
}

\begin{document}

\maketitle

\begin{abstract}
  Spectral projectors of second order differential operators play
  an important role in quantum physics and other scientific and
  engineering applications. In order to resolve local features and to
  obtain converged results, typically the number of degrees of freedom
  needed is much larger than the rank of the spectral projector. This
  leads to significant cost in terms of both computation and
  storage. In this paper, we develop a method to construct a basis
  set that is adaptive to the given differential operator. The basis
  set is systematically improvable, and the local features of the
  projector is built into the basis set. As a result the required
  number of degrees of freedom is only a small constant times the
  rank of the projector. The construction of the basis set uses a
  randomized procedure, and only requires applying the differential
  operator to a small number of vectors on the global domain, while
  each basis function itself is supported on strictly local domains and
  is discontinuous across the global domain. The spectral projector
  on the global domain is systematically approximated from such a
  basis set using the discontinuous Galerkin (DG) method. The global
  construction procedure is very flexible, and allows a local basis
  set to be consistently constructed even if the operator contains
  a nonlocal potential term.  We verify the effectiveness of the
  globally constructed adaptive local basis set using one-, two- and
  three-dimensional linear problems with local potentials, as well
  as a one dimensional nonlinear problem with nonlocal potentials
  resembling the Hartree-Fock problem in quantum physics.
\end{abstract}


\pagestyle{myheadings}
\thispagestyle{plain}

\section{Introduction}\label{sec:intro}

Consider the second order differential operator
\begin{equation}
  H = -\Delta + V(\vpos), \quad \vpos\in \Omega,
  \label{eqn:Hoperator}
\end{equation}
where $\Omega$ is a rectangular, bounded domain in $\RR^{d}$ with
periodic boundary conditions.  $V$ is a real, bounded and smooth
potential function. Then $H$ is a self-adjoint operator on $\Omega$.
Using the eigen-decomposition $H \psi_{i} = \varepsilon_{i} \psi_{i}$,
a spectral projector $P$ is an integral operator with its kernel
defined as
\begin{equation}
  P(\vpos,\vpos') = \mathbf{1}_{\mc{I}}(H)(\vpos,\vpos') =
  \sum_{\varepsilon_{i}\in \mc{I}} \psi_{i}(\vpos)\psi_{i}^{*}(\vpos').
  \label{eqn:projector}
\end{equation}
Here $\mc{I}$ is an interval that can be interpreted as an
energy window indicating the eigenfunctions of interest,
$\mathbf{1}_{\mc{I}}(\cdot)$ is an indicator function,
and $\psi_{i}^{*}(\vpos)$ is the complex conjugation of
$\psi_{i}(\vpos)$.  Denote by $n$ the number of eigenfunctions in
the summation of Eq.~\eqref{eqn:projector}, then the rank of $P$
is $n$. We assume that $n$ is large, which can range from hundreds
to hundreds of thousands. The spectral projector of such a form
or of similar forms arises in many scientific and engineering
problems. One notable example is the widely used Kohn-Sham density
functional theory~\cite{HohenbergKohn1964,KohnSham1965} in quantum
physics, where $\mc{I}$ contains the lowest $n$ eigenvalues of $H$.
Typically, a large number of degrees of freedom associated with a fine
numerical discretization are required to resolve the local features of
$\psi_{i}$'s with sufficient accuracy. This is the case when standard
discretization methods such as finite difference, finite element,
discontinuous Galerkin, planewave, and wavelet type of methods
are used. The ratio between the total number of degrees of freedom
(DOFs) and $n$ can range from hundreds to hundreds of thousands in
quantum physics applications~\cite{AlemanyJainKronikChelikowsky2004,
PaskKleinFongEtAl1999, GenoveseNeelovGoedeckerEtAl2008}.  As a
consequence, both the storage cost and the computation cost associated
with the spectral projector $P$ can be large.

\subsection{Contribution}
 
In this paper, instead of using a general basis set, we introduce
a new basis set that can be specifically tailored to represent the
spectral projector $P$, for a given operator $H$ and an interval
$\mc{I}$. The key observation is as follows. Let us partition $\Omega$
into a suitable collection of non-overlapping sub-domains called
\textit{elements}. If the size of each element is small enough,
then the numerical rank (a.k.a.  the approximate rank up to certain
truncation tolerance $\epsilon$~\cite{GolubVan2013}) of each row block
of $P$ restricted to any element can be bounded by a small constant.
We shall quantify the details of the statement above later in
the paper.  The singular value decomposition of one such row block
of $P$ defines the \textit{optimal} basis set on an element. Since
the local features of the range of $P$ is directly built into the
basis set, we can expect that the number of degrees of freedom in
such an optimal basis set is much smaller than that in a general
basis set.  However, such an optimal basis set cannot be practically
obtained, since it requires the knowledge of $P$ \textit{a priori}.
We devise a numerical algorithm to compute a \textit{nearly optimal}
basis set.  This is done by applying an approximate spectral projector,
characterized by a matrix function $f(H)$, to a small number of random
vectors defined on the global domain $\Omega$. The number of random
vectors is only slightly larger than the approximate rank of $f(H)$
restricted to each element.  The range of $P$ is then approximately
a subspace of the span of this basis set, and we find that this
is an efficient and accurate way to generate the basis functions
on all elements. Due to the non-overlapping condition, each basis
function is only supported on one element, and is discontinuous
on the global domain $\Omega$.  We use the discontinuous Galerkin
(DG) method~\cite{Arnold1982} to patch the basis set to obtain an
approximation to $\{\psi_{i}\}_{\varepsilon_{i}\in\mc{I}}$ or $P$.
Motivated by our previous work of the locally constructed adaptive
local basis set (LC-ALB) ~\cite{LinLuYingE2012, ZhangLinHuEtAl2017,
HuLinYang2015a},  the basis set in this work is dubbed the globally
constructed adaptive local basis set (GC-ALB). The LC-ALB
approach has already been demonstrated to be able to be executed
efficiently on massively parallel computers with over $100,000$
processes, and to efficiently perform large scale Kohn-Sham
density functional theory calculations for systems over $20,000$
atoms~\cite{HuLinYang2015a,BanerjeeLinHuEtAl2016}.

The GC-ALB set has the following advantages: 1) Systematically
improvable. As the number of basis functions in each element
increases, the accuracy of the projector represented in this basis
set systematically improves towards the converged spectral projector.
2) Efficient. The number of basis functions is directly related to the
numerical rank of the row blocks of the projector, and is much smaller
compared to the number of degrees of freedom needed to resolve the
local shape of $\{\psi_{i}\}$ in the real space.  The strict locality
of the basis set can significantly reduce the computation and storage
cost for $\{\psi_{i}\}$ and $P$.  3) Flexible. The construction
of the basis set only requires matrix-vector multiplication of
$H$ defined on the global domain $\Omega$. This allows existing
matrix-vector multiplication routines for computing $Hv$ to be readily
used without the need of constructing auxiliary operators. This
also facilitates generalizations to operators beyond the form
in~\eqref{eqn:Hoperator}. This can occur e.g.  for $H=-\Delta+V+W$,
where $W$ is an integral operator and hence $H$ becomes a nonlocal
operator. Such an operator arises in applications such as the
Kohn-Sham density functional theory with hybrid exchange-correlation
functionals~\cite{Becke1993,HeydScuseriaErnzerhof2003} and the
Hartree-Fock theory in quantum physics.

\subsection{Related work}

In the context of quantum physics, many tailored basis set have
been designed to reduce the number of DOFs to represent spectral
projectors (or density matrices in physics terminology).  Notable
examples include the Gaussian basis set and the atomic orbital basis
set~\cite{SolerArtachoGaleEtAl2002,Ozaki:03,BlumGehrkeHankeEtAl2009}.
Such basis sets are developed based on physical intuition, and
can provide relatively accurate solution with much reduced number
of degrees of freedom compared to more conventional basis sets.
However, expert knowledge is often required to systematically
converge the solution.  These basis sets have also been used
to ``enrich'' conventional basis sets to achieve a balance
between the small number of DOFs and the systematic convergence
property~\cite{SchwarzBlahaMadsen2002,SukumarPask2009}. However, the
number of DOFs in the mixed basis representation is often much larger
than those using Gaussian orbital or atomic orbital basis sets alone.

In order to achieve systematic convergence without sacrificing
the number of DOFs, one may give up the concept of designing
a basis set~\textit{a priori}, but instead generate the basis
set~\textit{on the fly}. This has been demonstrated via a number of
approaches based on filtration~\cite{LinLuYingE2012, LinYang2013,
RaysonBriddon2009, Garc'ia-CerveraLuXuanEtAl2009,MotamarriGavini2014}
as well as optimization~\cite{SkylarisHaynesMostofiEtAl2005,
LinLuYingE2012a,FattebertBernholc2000,MohrRatcliffBoulangerEtAl2014}
principles.  A common ingredient of these methods is to truncate the
$H$ operator into a series of operators defined only on different
sub-domains, and the basis set is then generated from the truncated
operators. This requires each basis function to satisfy zero Dirichlet
boundary conditions at each subdomain, which is not always achievable
without sacrificing the accuracy of the resulting basis set.
The method in this paper only uses matrix-vector multiplication
on the global domain, and hence the concern from the choice of
boundary conditions on local domains is completely removed. Our
numerical results indicate that the GC-ALB set can also be more
efficient than the LC-ALB set measured in terms of the number of
DOFs to reach the same target accuracy. The partition of unity
method (PUM)~\cite{BabuvskaMelenk1997} is another commonly used
option for enriching the basis set using local basis functions.
Compared to PUM, the basis functions in the DG approach are strictly
localized in non-overlapping elements, and hence are often better
conditioned. In fact, one can easily obtain an orthonormal basis set
in the DG approach through a local orthonormalization procedure. This
also facilitates the usage of efficient numerical techniques such
as the Chebyshev filtering techniques~\cite{ZhouSaadTiagoEtAl2006}
for the computation of spectral projectors.

\subsection{Outline of the paper}

The rest of the paper is organized as follows. We review the interior
penalty formulation of the discontinuous Galerkin framework,
introduce the optimal discontinuous basis set and the locally
constructed adaptive local basis set in section~\ref{sec:prelim}. We
present the globally constructed adaptive local basis set in
section~\ref{sec:gcalb}. The numerical results are given in
section~\ref{sec:numer}, followed by the conclusion and discussion
in section~\ref{sec:conclusion}.

\section{Preliminaries}\label{sec:prelim}

\subsection{Discontinuous Galerkin method}\label{sec:dg}

Without loss of generality, let $\Omega=(0,L)^d$ where
$d=1,2,3$, and $\mc{K}$ be a regular partition of $\Omega$ into a set
of non-overlapping elements.  For $\kappa\in \mc{K}$, we denote by
$\overline{\kappa}$ the closure of $\kappa$.  For any two elements
$\kappa,\kappa'\in \mc{K}$, The periodic boundary condition on
$\Omega$ implies that the partition is regular across the boundary
$\partial \Omega$.  We remark that generalizations to other boundary
conditions such as Dirichlet or Neumann boundary conditions, as well
as to non-rectangular domains, can be done with minor modification.

We denote by $H^1(\kappa)$ the standard Sobolev space of
$L^2(\kappa)$-functions such that the first partial derivatives are
also in $L^2(\kappa)$.  We denote the set of piecewise $H^1$-functions
by
\[
	H^1(\mc{K}) = \left\{ v\in L^2(\Omega) \,\middle|\, v|_\kappa \in
  H^1(\kappa), \quad \forall \kappa \in \mc{K} \right\},
\]
which is also referred to as the broken Sobolev space.  For $v,w\in
H^{1}(\mc{K})$, the inner product is
\begin{equation}
  (v,w)_{\mc{K}} = \sum_{\kappa\in\mc{K}} (v,w)_\kappa :=
  \sum_{\kappa\in\mc{K}} \int_{\kappa} v^{*}(\vpos) w(\vpos) \ud \vpos,
  \label{eqn:innerprod}
\end{equation}
which induces a norm $\|v\|_\mc{K} = (v,v)_\mc{K}^\frac12$.

For $v,w\in H^{1}(\mc{K})$ and $\kappa,\kappa' \in \mc{K}$,
define the jump and average operators
on a face $\overline{\kappa}\cap\overline{\kappa}'$ by
\begin{equation}
  \avg{v} = \tfrac12(v|_{\kappa} + v|_{\kappa'}), \quad 
  \avg{\nabla v} = \tfrac12(\nabla v|_{\kappa} + \nabla v|_{\kappa'}),
  \label{}
\end{equation}
and
\begin{equation}
  \jmp{v} = v|_{\kappa} \vn_\kappa + v|_{\kappa'} \vn_{\kappa'}, \quad
  \jmp{\nabla v} = \nabla v|_{\kappa} \cdot \vn_\kappa +\nabla
      v|_{\kappa'} \cdot \vn_{\kappa'},
  \label{}
\end{equation}
where $\vn_\kappa$ denotes the exterior unit normal of the element
$\kappa$.

In order to numerically solve the eigenvalue problem 
\[
H \psi_{i}=\varepsilon_{i} \psi_{i},
\]
we need to identify a basis set which spans a subspace
of $H^{1}(\mc{K})$.  Let $N_\kappa$ be the number
of DOFs on $\kappa$, and the total number of DOFs
is $N_{\mc{K}}=\sum_{\kappa\in\mc{K}} N_{\kappa}$.  Let
$\VN(\kappa)=\text{span}~\{\varphi_{\kappa,j}\}_{j=1}^{N_{\kappa}}$,
where each $\varphi_{\kappa,j}$ is a function defined on $\Omega$
with compact support only in $\kappa$. Hence $\VN(\kappa)$
is a subspace of $H^{1}(\mc{K})$ and is associated with a
finite dimensional approximation for $H^{1}(\kappa)$. Then
$\VN=\bigoplus_{\kappa\in\mc{K}}\VN(\kappa)$ is a finite dimensional
approximation to $H^{1}(\mc{K})$. We also assume all functions
$\{\varphi_{\kappa,j}\}$ form an orthonormal set of vectors in the
sense that
\begin{equation}
  (\varphi_{\kappa,j},\varphi_{\kappa',j'})_{\mc{K}} =
  \delta_{\kappa,\kappa'}\delta_{j,j'}, \quad \forall \kappa,\kappa'\in
  \mc{K},  1\le j \le N_{\kappa},  1\le j' \le N_{\kappa'}.
  \label{eqn:orthonormal}
\end{equation}

The interior penalty formulation of the discontinuous Galerkin
method~\cite{Arnold1982} introduces the following bilinear form
\begin{multline}
  \label{eqn:bilinear}
	a(w,v) = 
		\sum_{\kappa\in\mc{K}}
		\Big[
				(\nabla w, \nabla v)_\kappa 
				+ (Vw,v)_\kappa 
		\Big]
		+
		\tfrac12		
		\sum_{\kappa\in\mc{K}}
		\Big[
				- (\nabla w,\jmp{v})_{\partial\kappa} 
				- (\jmp{w},\nabla v)_{\partial\kappa} \Big]\\
    + \tfrac12		
		\sum_{\kappa\in\mc{K}} \Big[
\gamma_{\kappa} ( \jmp{w},\jmp{v})_{\partial\kappa}\Big].
\end{multline}
Here the terms in the first bracket corresponds to the
operator $H$. The terms in the second bracket are obtained
from integration by parts, and the terms in the third bracket
is a penalty term to guarantee the stability of the bilinear
form~\cite{ArnoldBrezziCockburnEtAl2002}.  The penalty parameter
$\gamma_{\kappa}$ on each element $\kappa$ needs to be large enough,
and the value of $\gamma_{\kappa}$ depends on the choice of basis
set $\VN$. For general non-polynomial basis functions the value
of $\gamma_{\kappa}$ is not known \textit{a priori}. One possible
solution is given in~\cite{LinStamm2016} which provides a formula for
evaluating $\gamma_{\kappa}$ on the fly for general non-polynomial
basis sets based on the solution of eigenvalue problems restricted
to each element $\kappa$.

Using the bilinear form~\eqref{eqn:bilinear},
the solution of 
\begin{equation}
	\label{eqn:DGeig1}
    a(\psi^{\VN}_{i},v) = \varepsilon^{\VN}_{i}
    (\psi^{\VN}_{i},v)_\mc{K},\qquad\forall v\in\VN
\end{equation}
gives the numerical solution of eigenpairs of the form
$(\varepsilon^{\VN}_{i},
\psi^{\VN}_{i})$ and $\psi^{\VN}_{i}\in H^{1}(\mc{K})$.
Eq.~\eqref{eqn:DGeig1} can be equivalently written as a standard
linear eigenvalue problem
\begin{equation}
  \sum_{\kappa',j'} H^{\VN}_{\kappa,j;\kappa',j'} c_{\kappa',j';i} =
    \varepsilon^{\VN}_{i} c_{\kappa,j;i},
  \label{eqn:DGeig}
\end{equation}
where $\{c_{\kappa,j;i}\}$ satisfies $\psi^{\VN}_{i} =
\sum_{\kappa,j}c_{\kappa,j;i}\varphi_{\kappa,j}$,
and the reduced matrix $H^{\VN}$ is of size $N_{\mc{K}}\times
N_{\mc{K}}$ with matrix
elements
\begin{equation}
  H^{\VN}_{\kappa,j;\kappa',j'}  =
  a(\varphi_{\kappa,j},\varphi_{\kappa',j'}).
  \label{eqn:Helement}
\end{equation}
Using the solution of Eq.~\eqref{eqn:DGeig}, we can select
$\varepsilon^{\VN}_{i}\in \mc{I}$ and obtain an approximation to the
spectral projector
\begin{equation}
    P(\vpos,\vpos') \approx \sum_{\varepsilon^{\VN}_{i} \in \mc{I}}
    \psi^{\VN}_{i}(\vpos) \left( \psi^{\VN}_{i}(\vpos')\right)^*
    = \sum_{\kappa,\kappa',j,j'}
  \varphi_{\kappa,j}(\vpos) \Gamma_{\kappa,j;\kappa',j'}
  \varphi^{*}_{\kappa',j'}(\vpos').
  \label{eqn:projectordlr}
\end{equation}
Here $\Gamma$ is the $N_{\mc{K}}\times N_{\mc{K}}$ matrix
representation of $P$ in the basis set $\VN$, and
\begin{equation}
    \Gamma_{\kappa,j;\kappa',j'} = \sum_{\varepsilon^{\VN}_{i}\in\mc{I}}
  c_{\kappa,j;i} c^{*}_{\kappa',j';i}.
  \label{eqn:expandP}
\end{equation}

\subsection{Optimal discontinuous basis set}\label{sec:Optbasis}

The discontinuous Galerkin method in section~\ref{sec:dg} can be
applied to very general basis sets $\VN$.  Here we consider the
\textit{optimal} basis set $\VN$ for representing the spectral
projector $P$ with a discontinuous basis set. To simplify the
discussion below, we also use linear algebra notation in this
section when necessary.  This means that we may not distinguish the
kernel of an operator and a finite dimensional matrix consisting
of its nodal values discretized on a fine set of real space grid
points, with the number of grid points denoted by $N_{g}$. Then
notation such as $\vpos,\vpos'$ can be real space grid points
in $\Omega$, or row / column indices of vectors / matrices. We
call $P(\vpos,:):=\{P(\vpos,\vpos'),\vpos'\in \Omega\}$
a row vector, and $P(:,\vpos):=\{P(\vpos',\vpos),\vpos'\in
\Omega\}$ a column vector, respectively. Similarly, we call
$P(\kappa,:):=\{P(\vpos,\vpos'),\vpos\in \kappa,\vpos'\in \Omega\}$
a row block, and $P(:,\kappa):=\{P(\vpos',\vpos),\vpos\in
\kappa,\vpos'\in \Omega\}$ a column block, respectively.

\begin{figure}[h]
  \begin{center}
    \includegraphics[width=0.9\textwidth]{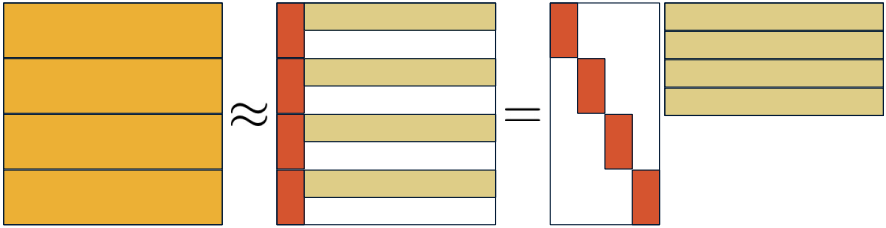}
  \end{center}
  \caption{Construction of the optimal discontinuous basis set for spectral
  projector. Left: the spectral projector $P$ is partitioned into 4 row
    blocks; Middle:
    each $P(\kappa,:)$ is low-rank factorized via SVD, i.e., $P(\kappa,:)
    \approx \Phi_{\kappa} (S_\kappa V^*_\kappa)$; Right: subspace $\VN$ is
    assembled from $\{\Phi_{\kappa}\}$.} \label{fig:BasisFact}
\end{figure}

Since the rank of the spectral projector $P$ is $n$, if we choose a
partition $\mc{K}$ fine enough we may expect that the numerical rank
of each row block $P(\kappa,:)$ becomes small.  In particular,
note that the rank of $P(\kappa,:)$ cannot exceed the number of degrees of
freedom in $\kappa$, which is a constant and is independent of $n$. Our
numerical results indicate that this rank can often be much lower
than the number of degrees of freedom in $\kappa$ in practice. The
singular value decomposition (SVD) of $P(\kappa,:)$ can be written as
\begin{equation}
  P(\kappa,:) \approx \Phi_{\kappa} S_{\kappa} V_{\kappa}^*,
  \label{eqn:optimalP}
\end{equation}
where $S_{\kappa}$ is a diagonal matrix containing the
leading $N_{\kappa}$ singular values on $\kappa$, and
$\Phi_{\kappa}(\vpos) = [\varphi_{\kappa,1}(\vpos), \ldots,
\varphi_{\kappa,N_{\kappa}}(\vpos)]$ for $\vpos \in \kappa$.
The support of each function $\varphi_{\kappa,j}\in H^{1}(\mc{K})$
is strictly in $\kappa$.  Since $\Phi_{\kappa}$'s are generated
from the SVD of $P$, clearly the range of $P$ is approximately
contained in $\text{span} \{\varphi_{\kappa,j}\}$. For a given
$\kappa\in \mc{K}$, the basis $\Phi_{\kappa}$ achieves the smallest
error in 2-norm for representing $P(\kappa,:)$ thanks to the
optimal approximation property of the SVD~\cite{GolubVan2013}
using $N_{\kappa}$ basis functions. Hence the basis set
$\{\Phi_{\kappa}\}_{\kappa\in\mc{K}}$ can be regarded as the optimal
discontinuous basis set for representing $P$ for a given set of degrees
of freedom $\{N_{\kappa}\}_{\kappa\in\mc{K}}$.  We illustrate the
decomposition~\eqref{eqn:optimalP} for the entire spectral projector
$P$ in Fig.~\ref{fig:BasisFact}.

\subsection{Locally constructed adaptive local basis set}\label{sec:lcalb}

The optimal discontinuous basis set cannot be used for practical
computation, since its construction depends on the knowledge
of $P$. One possible approximation of such a basis set using
non-polynomial basis functions is the adaptive local basis (ALB)
set~\cite{LinLuYingE2012}.  More specifically, we refer to this
basis set the locally constructed adaptive local basis (LC-ALB) set,
in order to distinguish from the globally constructed adaptive local
basis set in section \ref{sec:gcalb}.

\begin{figure}[h]
  \begin{center}
    \includegraphics[width=0.37\textwidth]{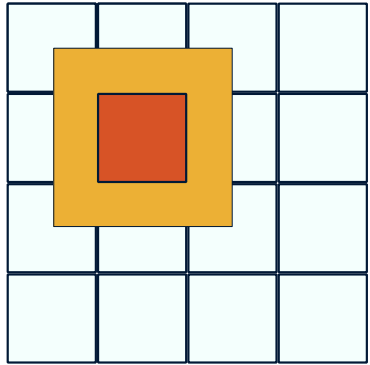}
  \end{center}
  \caption{The entire two-dimensional domain is partitioned into 4 by 4
    blocks denoted by the white blocks. A particular element $\kappa$ and
    its buffer region $\widetilde{\kappa} \backslash \kappa$ are denoted as
    the red block and yellow block respectively.}
  \label{fig:LCALBDomain}
\end{figure}

Consider the case that $\mc{I}$ contains the lowest $n$ eigenvalues
of $H$. In the $d$-dimensional space, for each element $\kappa$,
we form an \textit{extended element} $\widetilde{\kappa}$ around
$\kappa$, and we refer to $\widetilde{\kappa}\backslash \kappa$
as the buffer region for $\kappa$. Fig.~\ref{fig:LCALBDomain}
illustrates a particular $\kappa$ together with its butter region.
On $\widetilde{\kappa}$ we solve the eigenvalue problem
\begin{equation}
  -\Delta \widetilde{\varphi}_{i} + V \widetilde{\varphi}_{i} =
  \lambda_{i} \widetilde{\varphi}_{i}, \label{eqn:localeig}
\end{equation}
with certain boundary conditions on $\partial \widetilde{\kappa}$.
This eigenvalue problem can be solved using standard basis set such as
finite difference, finite elements, or planewaves.  For the numerical
examples in this paper, the periodic boundary conditions is applied
on each $\partial \widetilde{\kappa}$, and the eigenvalue problem
is solved via the pseudo-spectral method (the planewave basis set).
Note that the size of the extended element $\widetilde{\kappa}$
is independent of the size of the global domain, and so is the
number of basis functions per element. In order to obtain $\VN$,
the eigenfunctions corresponding to lowest $N_{\kappa}$ eigenvalues
are restricted from $\widetilde{\kappa}$ to $\kappa$, i.e.
\[
\varphi_{i}(x) = \begin{cases}
  \left[\widetilde{\varphi}_{i}\right]\vert_{\kappa}(x), & x\in \kappa;\\
  0,&\text{otherwise},
\end{cases} \quad i=1,\ldots,N_{\kappa}.
\]
After orthonormalizing $\{\varphi_{i}\}$ locally on each element
$\kappa$ and removing the linearly dependent functions via a local
singular value decomposition, the resulting set of orthonormal
functions form the LC-ALB set.

The advantage of the LC-ALB set is that the basis functions for each
element $\kappa$ can be generated completely independently. However,
due to the fictitious boundary conditions imposed on the extended
element $\partial \widetilde{\kappa}$, the effectiveness of the
LC-ALB set depends on the size of the buffer region. On one extreme,
if the size of the buffer region is $0$ and when the periodic boundary
condition is used, since $V$ is in general not a periodic function
on $\kappa$, the accuracy of the basis set can be severely affected
by the Gibbs phenomena. On the other hand, if the buffer region is
chosen to be too large, then the solution of the local eigenvalue
problem~\eqref{eqn:localeig} can become expensive. In practice we
find that choosing $\widetilde{\kappa}$ to contain $\kappa$ and
its $3^{d}-1$ neighboring elements yields a relatively good balance
between efficiency and accuracy, as has been demonstrated by the usage
for solving PDEs~\cite{LinStamm2016,LinStamm2017} and for solving
practical Kohn-Sham equations~\cite{LinLuYingE2012, HuLinYang2015a}.


\section{Globally constructed adaptive local basis set}\label{sec:gcalb}

In this section, we propose a new strategy to construct an
approximation to the optimal discontinuous basis set by using
matrix-vector multiplication involving the operator $H$ defined on the
global domain $\Omega$. This allows us to overcome the difficulty of
choosing the boundary condition and the size of the extended element
as in the LC-ALB set. Numerical results indicate that the resulting
basis set is more effective in terms of the number of DOFs, and the
strategy can be adapted to more general cases such as when $V$ is an
integral operator with a nonlocal kernel.

\subsection{Formulation}

We first introduce Algorithm~\ref{alg:randomsvd}, which is a variant
of e.g. Algorithm 4.1 in~\cite{HalkoMartinssonTropp2011} for finding
the approximate range of a numerically low rank matrix.

\begin{algorithm}[H]
  \caption{Randomized range finder for a given matrix $A$.}

  Input: $A\in \CC^{p\times q}$. Approximate rank $k$.

  Output: Left-singular vectors $U\in \CC^{p\times k}$.

\begin{algorithmic}[1]

  \STATE Generate an orthonormal random matrix $R\in
  \CC^{q\times (k+c)}$ where $c$ is a small oversampling constant.

  \STATE Compute $W=A R$. \label{step:matvec}

  \STATE Perform the SVD for $W=USV^{*}$, with the
  diagonal entries of $S$ ordered non-increasingly.
  
  \STATE Return the first $k$ columns of $U$.
  
\end{algorithmic}
\label{alg:randomsvd}
\end{algorithm}

If we treat $A$ as a dense matrix and apply the SVD directly, the
computational complexity will be $\Or(pqk)$. On the other hand,
Alg.~\ref{alg:randomsvd} only requires applying the matrix $A$
to $(k+c)$ vectors, together with the SVD for $W$ which costs
$\Or(pk^2)$ operations. Hence the randomized range finder algorithm
significantly reduces the cost, if $k$ is much smaller than $q$
and if the matrix vector multiplication $Av$ can be evaluated
quickly. Usually, step~\ref{step:matvec} is the most expensive
operation in Algorithm~\ref{alg:randomsvd}.

Assume $\mc{K}$ is a partition of $\Omega$ so that each matrix row
block $P(\kappa,:)$ is a numerically low rank matrix.  If we apply
Algorithm~\ref{alg:randomsvd} to $P(\kappa,:)$, the output gives highly
accurate approximation to the optimal basis set $\{\Phi_{\kappa}\}$
for $\kappa\in \mc{K}$.  Furthermore, the random matrix $R$ can
be repeatedly used for different $\kappa \in \mc{K}$. Therefore,
the matrix-vector multiplication for different matrix row blocks
$P(\kappa,:)$ do not need to be applied independently. Instead it
is equivalent to apply the entire matrix $P$ to a random matrix
$R$, and to perform the SVD for each element independently to
obtain an approximate range represented by $\Phi_{\kappa}$ for each
$P(\kappa,:)$. The collection of the functions $\{\Phi_{\kappa}\}$
gives the globally constructed adaptive local basis set (GC-ALB).
Algorithm~\ref{alg:gcalb} describes this procedure for a general matrix
$A\in \CC^{N_{g}\times N_{g}}$, where $N_{g}$ is the number of DOFs
corresponding to a fine numerical discretization such as planewaves.

\begin{algorithm}[h]
  \caption{Globally constructed adaptive local basis set for a given
  matrix $A$.}

  Input: $A\in \CC^{N_{g}\times N_{g}}$. Partition $\mc{K}=\{\kappa\}$
  with approximate rank for each element $\{N_{\kappa}\}$.

  Output: The basis set $\{\Phi_{\kappa}\}$.

\begin{algorithmic}[1]

  \STATE Generate an orthogonal random matrix $R\in
  \CC^{N_{g}\times (\max_{\kappa}N_{\kappa}+c)}$, where $c$ is a small
  oversampling constant.

  \STATE Compute $W=A R$.\label{step:matvec_gcalb}

  \FOR {$\kappa\in \mc{K}$}
  \STATE Perform the SVD for
  $W(\kappa,:)=U_{\kappa}S_{\kappa}V_{\kappa}^{*}$, with the
  diagonal entries of $S_{\kappa}$ ordered
  non-increasingly.\label{step:svd_gcalb}
  
    \STATE Obtain $\Phi_{\kappa}$ from the first $N_{\kappa}$ columns of
  $U_{\kappa}$.
  \ENDFOR
  
\end{algorithmic}
\label{alg:gcalb}
\end{algorithm}

\begin{figure}[h]
  \begin{center}
    \includegraphics[width=0.9\textwidth]{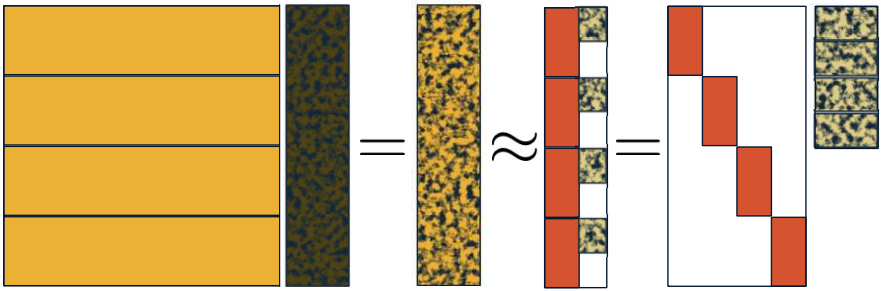}
  \end{center}
  \caption{Globally constructed adaptive local basis functions. Left:
  the spectral
    projector $P$ is applied to random vectors $R$; Middle-left: the
    result of the left part contains the column basis of each block;
    Middle-right: the column basis of each block is revealed via SVD;
    Right: subspace $\VN$ is assembled by $\{\Phi_\kappa\}$.}
  \label{fig:BasisFactRand}
\end{figure}

When taking the matrix $A$ to be the spectral projector $P$,
Fig.~\ref{fig:BasisFactRand} illustrates Alg.~\ref{alg:gcalb} for
the case that $\mc{K}$ is partitions the domain $\Omega$ into $4$
elements.  Comparing to Fig.~\ref{fig:BasisFact} where each block
of $P$ is explicitly factorized, Fig.~\ref{fig:BasisFactRand} first
applies $P$ to random vectors $R$ and then factorizes each block of
$PR$. Such an extra step is crucial here. Computing the dense $P$
is expensive in terms of both computation and memory, whereas the
matrix vector multiplication $Pv$ can be efficiently calculated,
which could be orders of magnitudes faster for large problems.

\subsection{Rational approximation for matrix-vector multiplication}
\label{sec:RationalApprox}

In order to construct the GC-ALB set for the projector $P$,
Alg.~\ref{alg:gcalb} requires an efficient method to compute the
matrix-vector multiplication $Pv$. Since the spectral projector
is a non-smooth matrix function $\mathbf{1}_{\mc{I}}(\cdot)$,
the computation of $Pv=\mathbf{1}_{\mc{I}}(H)v$ may still be a
costly procedure. Fortunately, we only need Alg.~\ref{alg:gcalb}
to find an approximate range of $P$. Hence we may replace
$\mathbf{1}_{\mc{I}}(\cdot)$ by a smooth function $f(\cdot)$,
with the requirement that the support of $f$ covers the interval
$\mc{I}$, and that $f(H)v$ is relatively easy to compute. Then we can
apply Alg.~\ref{alg:gcalb} to find the approximate range of $f(H)$.
The choice of $f$ is certainly not unique. Here we use a modified
Zolotarev's function to be $f(\cdot)$, which is an optimal rational
approximation to the indicator function as to be demonstrated below.

Without loss of generality, we assume that $\mc{I}=[a,b]$ in the
following discussion.  Zolotarev's function $Z_{2r}(x;\ell)$ was
initially proposed as the best rational approximant of type $(2r-1,2r)$
for the signum function $\sign(x)$ on the interval $[-1,-\ell] \cup
[\ell,1]$ with a positive parameter $\ell<1$~\cite{Zolotarev1877,
Akhiezer1990}. Recently, it was composed with a M\"obius transformation
$T(x)$ and a linear transformation~\cite{Guttel2015a, Li2017b} to
become the best rational approximant of type $(2r,2r)$ for an indicator
function $\mathbf{1}_{[a,b]}(x)$ on the interval $(-\infty,a_-] \cup
[a,b] \cup [b_+,+\infty)$, where $a_-$ and $b_+$ are two parameters
such that $a_- < a$ and $b < b_+$.  Both the M\"obius transformation
$T(x)$ and the parameter $\ell$ in Zolotarev's function depend on
$a_-$ and $b_+$. To be more precise, the M\"obius transformation is
defined as follows,
\begin{equation} \label{eqn:Mobius}
    T(x) = \gamma \frac{ x - \alpha }{ x - \beta },
\end{equation}
with $\alpha \in (a_-,a)$ and $\beta \in (b,b_+)$ such that
\begin{equation} \label{eqn:MobiusEqs}
    T(a_-) = -1,\quad T(a) = 1, \quad T(b) = \ell, \text{ and } T(b_+) =
    -\ell.
\end{equation}
Here, the variables $\alpha, \beta, \gamma,$ and $\ell$ are
determined by $a_-, a, b,$ and $b_+$ via solving the equations in
\eqref{eqn:MobiusEqs}. Combining with a simple linear transformation,
$(x+1)/2$, we arrive at a modified Zolotarev's function,
\begin{equation} \label{eqn:MZoloFunc}
    \begin{split}
	R(x) = & \frac{Z_{2r}(T(x);\ell) + 1}{2}\\
    = &
	\frac{M}{2}\sum_{j=1}^r \frac{a_j \gamma}{\gamma^2 + c_{2j-1}} +
	\frac{1}{2} + \frac{M}{2}\sum_{j=1}^r \left(\frac{w_j}{x-\sigma_j} +
	\frac{\bar{w_j}}{x - \bar{\sigma}_j}\right),
    \end{split}
\end{equation}
where $\gamma$ is the same as in \eqref{eqn:Mobius}, $M, a_j, c_j,
w_j,$ and $\sigma_j$ are constants as defined in~\cite{Li2017b},
and $\bar{\cdot}$ denotes the complex conjugate. $\left\{\sigma_j,
\bar{\sigma}_j\right\}_{j=1}^r$ are known as the poles of the modified
Zolotarev's function.

When the modified Zolotarev's function is used as $f(\cdot)$, and the
matrix $A$ in Alg.~\ref{alg:gcalb} is replaced by $f(H)$, the line $2$
in Alg.~\ref{alg:gcalb} can be evaluated via,
\begin{equation}\label{eqn:MZoloFuncMat}
    \begin{split}
        f(H)R = & \left( \frac{M}{2}\sum_{j=1}^r
        \frac{a_j \gamma}{\gamma^2 + c_{2j-1}} + \frac{1}{2} \right) R \\
        & + \frac{M}{2}\sum_{j=1}^r
        \left(w_j\left(H - \sigma_j I \right)^{-1} R + \bar{w_j}\left(H -
        \bar{\sigma}_j I \right)^{-1} R \right).
    \end{split}
\end{equation}
This requires solving $2r$ complex-shifted linear systems, where
$I$ denotes the identity matrix of the same size as $H$. If
both $H$ and $R$ are real matrices, then $\bar{w}_j\left(H -
\bar{\sigma}_j I \right)^{-1}R$ is the complex conjugate of
$w_j\left(H-\sigma_j\right)^{-1}R$. Therefore, solving $2r$
shifted linear systems in \eqref{eqn:MZoloFuncMat} can be reduced
to solving $r$ shifted linear systems instead.  These shifted
linear systems can be solved via standard iterative methods such as
GMRES~\cite{SaadSchultz1986} and MINRES~\cite{PaigeSaunders1975} with
a preconditioner.  The condition number of the shifted linear
systems depends on the minimal distance between the eigenvalues of
$H$ and the shifts $\{\sigma_{j},\bar{\sigma}_j\}$, as well as the
spectral radius of $H$.  The minimal distance can be systematically
controlled by tuning the smoothness of the function $f(\cdot)$,
through the adjustment of the interval $(a_-,a)$ and $(b,b_+)$.
When the spectral radius of $H$ is large such as in the case of
the planewave discretization, we find that the inverse of a shifted
Laplacian~\cite{GijzenErlanggaVuik2007} is an efficient preconditioner
to reduce the condition number.  Therefore, as shown in the numerical
results, the number of iterations needed for solving the linear
systems can be systematically controlled and relatively small.

\begin{figure}[h]
  \begin{center}
    \subfloat[]{\includegraphics[width=0.44\textwidth]{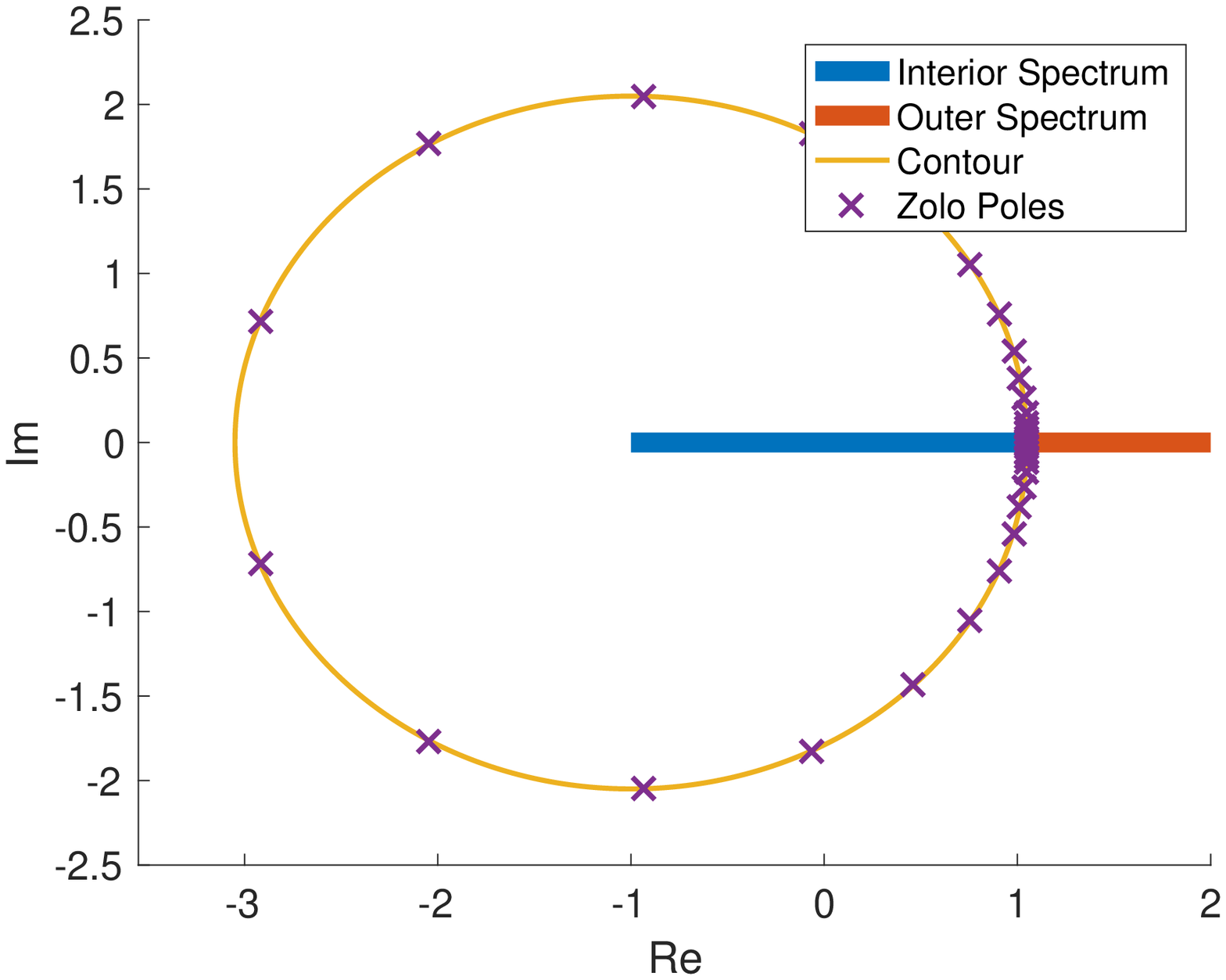}}~~
    \subfloat[]{\includegraphics[width=0.45\textwidth]{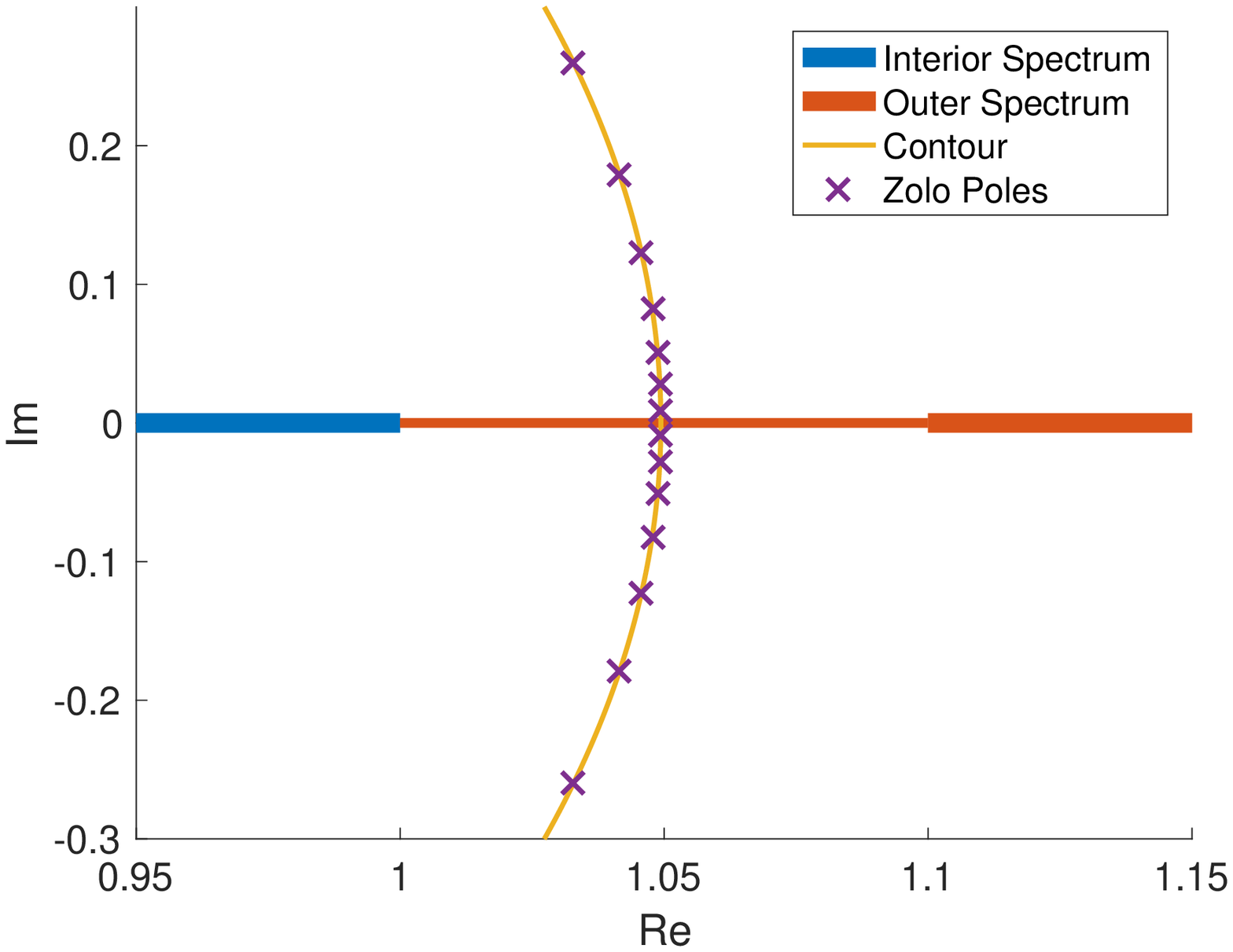}}\\
    \subfloat[]{\includegraphics[width=0.45\textwidth]{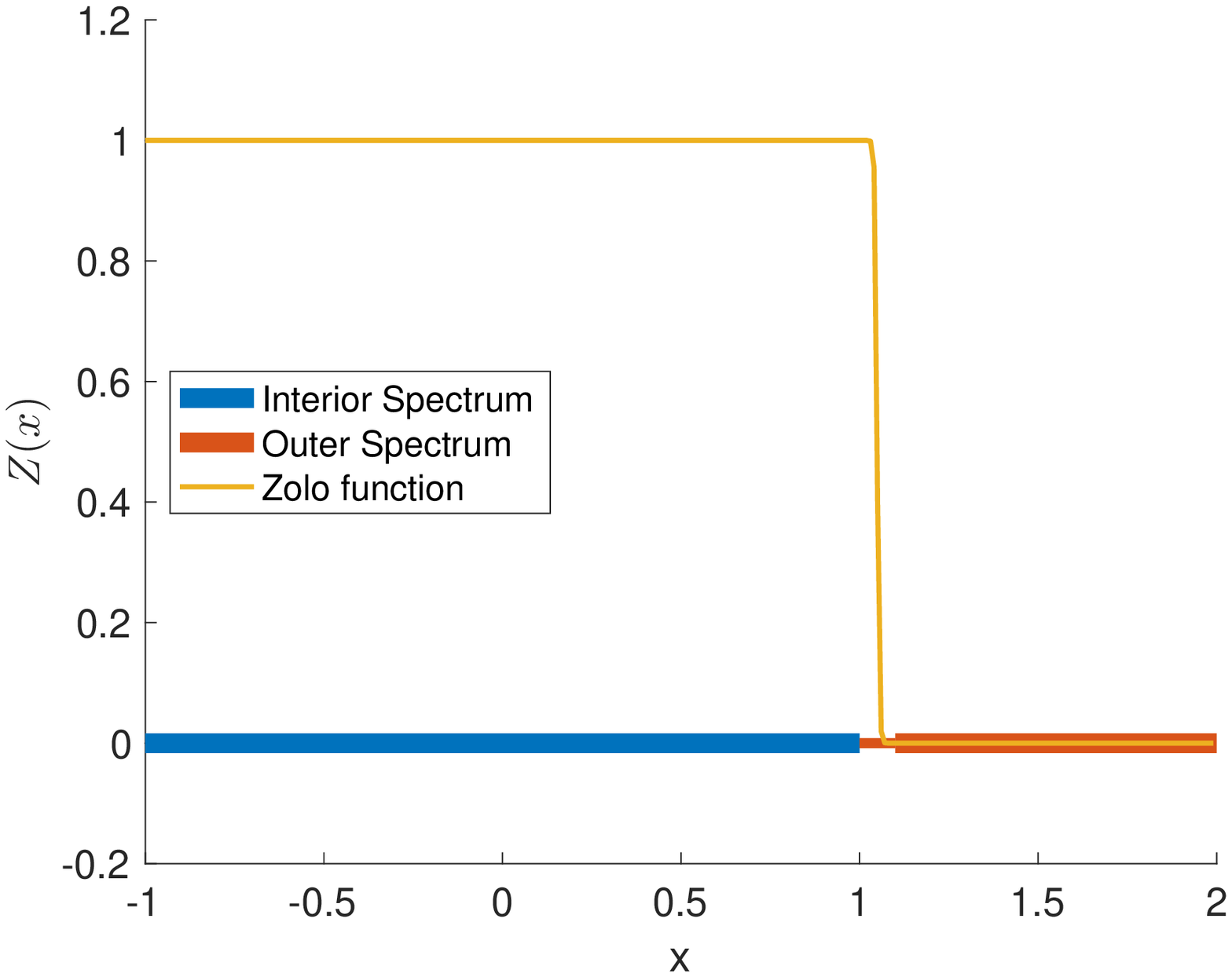}}~~
    \subfloat[]{\includegraphics[width=0.50\textwidth]{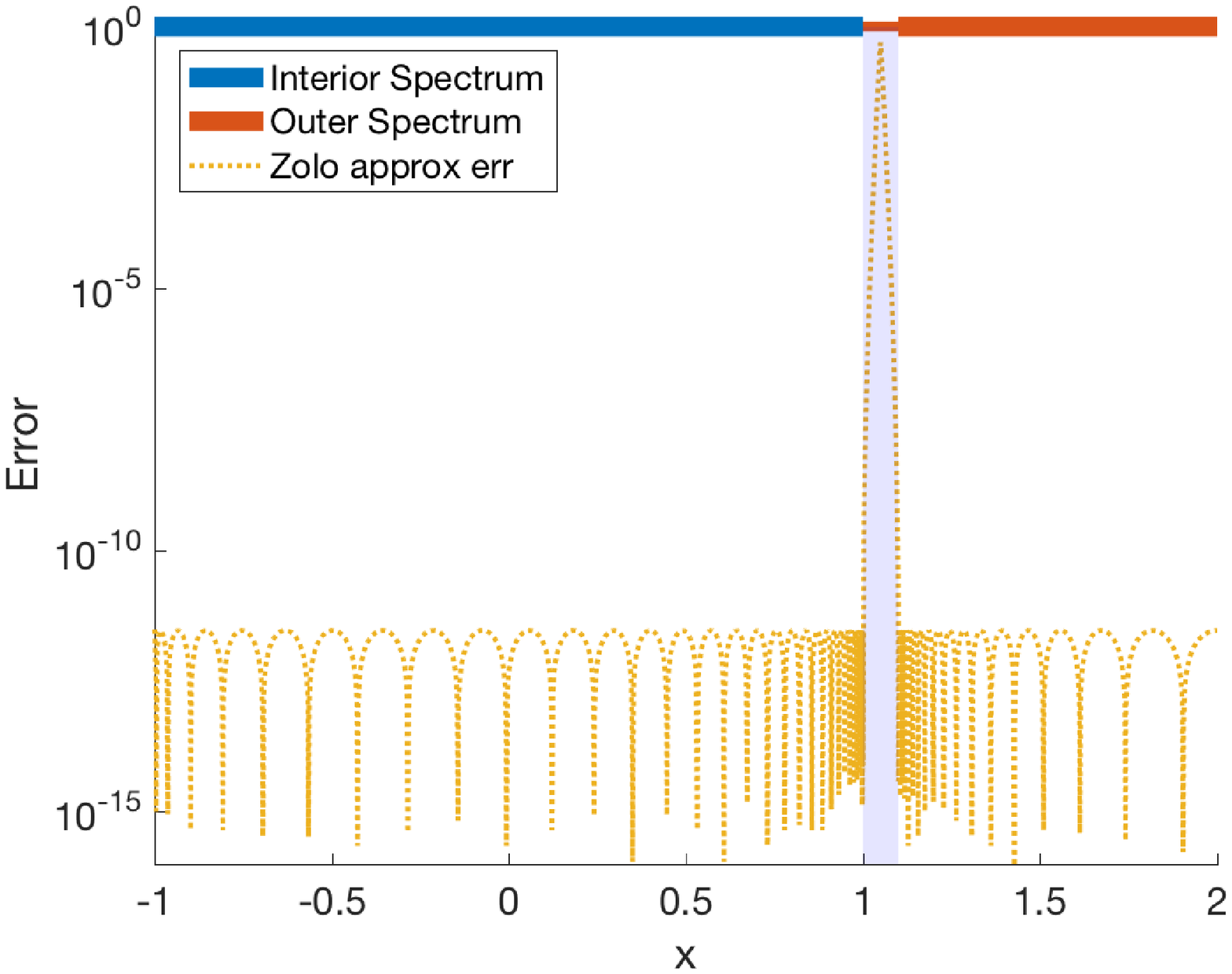}}
  \end{center}
    \caption{An example with $\mc{I} = [-1,1]$ and $(1,\infty)$ being the
    excluded spectrum. A modified Zolotarev's function is constructed with
    an artificial gap $(1,1.1)$ and $r=16$. (a) and (b) describe the
    locations of the poles on the complex plane. (c) is the modified
    Zolotarev's function and (d) is the absolute error comparing to
    indicator function $\mathbf{1}_{\mc{I}}(\cdot)$.}
  \label{fig:ZoloApprox}
\end{figure}

Through the discussion above, the choice of $a_-$ and $b_+$ remains to
be determined. For a fixed indicator function $\mathbf{1}_{\mc{I}}$
with the given interval $\mc{I} = [a,b]$, $a_-$ and $b_+$ determine
the quality of the approximation of the modified Zolotarev's function
in \eqref{eqn:MZoloFunc}.  Generally, if either interval $[a_-,a]$ or
$[b,b_+]$ becomes too narrow, it may require a large number of terms
$r$ in \eqref{eqn:MZoloFuncMat} to reach the same target accuracy. This
translates to solving more shifted linear systems. However,
the situation simplifies when $\mc{I}=[a,b]$ covers the lowest
$n$ eigenvalues of $H$, as will be demonstrated in the numerical
results. Let the eigenvalues of $H$ be $\lambda_1 \leq \lambda_2 \leq
\dots \leq \lambda_n \leq \lambda_{n+1} \leq \dots \leq \lambda_{N_g}$.
The parameter $a_-$ can be an arbitrary number in $(-\infty, a)$,
and we can set $a_-$ to be $-\infty$. The choice of $b_+$ relies on
the spectrum property of $H$ around $b$. When there is an eigenvalue
gap around $b$, i.e., $\lambda_n \leq b < \lambda_{n+1}$, $b_+$ is
set to be $\lambda_{n+1}$ or its estimated lower bound calculated
via a few steps of Lanczos method~\cite{ZhouSaadTiagoEtAl2006}. In
the case that $H$ has continuous spectrum around $b$, we construct a
small gap as $b_-=b$ and $b_+ = b+\delta$ for some small positive
constant $\delta$.  The consequence of such a gap is that the
approximated projector of $P$ would include extra eigenvectors with
non-zero weights. In practice, we find that the GC-ALB method is
robust to this choice of $\delta$.  Such an observation even allows
us to choose $b_{+}$ to be larger than $\lambda_{n+1}$ even in the
presence of a gap, in order to reduce $r$ and hence the computational
cost.  When the location of $\lambda_{1},\lambda_{n}$ is not known
\textit{a priori}, similar to the situation in Chebyshev filtering
techniques~\cite{ZhouSaadTiagoEtAl2006,ZhouChelikowskySaad2014},
the initial guess of $(a,b)$ can be efficiently obtained through a
few Lanczos~\cite{Lanczos1950} iterations in practice.

Fig.~\ref{fig:ZoloApprox} gives an example of the modified
Zolotarev's function for the approximation of $\mathbf{1}_{\mc{I}}$
where $\mc{I}=[-1,1]$. We assume $1$ is in the continuous spectrum
of $H$. We choose $a_- = -\infty$ and $b_{+} = 1.1$ so that $r=16$
is sufficient to approximate the indicator function with error below
$10^{-10}$ in the interval $(a,b) \cup (b_+,\infty)$.

\subsection{Complexity}

In practical computation of the spectral projector, the following
two scenarios are often encountered when counting the complexity with
respect to the increase of the number of DOFs.
\begin{enumerate}
  \item The size of the global domain $\Omega$ is fixed, and the
  number of DOFs increases due to the refinement of the discretization.
  \item The size of the global domain $\Omega$ increases, and the
  number of DOFs increases proportionally to the volume of $\Omega$.
\end{enumerate}

Let $M$ be the number of elements in $\mc{K}$, and $N_{g}$ be the
number of DOFs corresponding to a fine discretization on the global
domain. For simplicity let all elements have the same number of basis
functions denoted by $n_b$, and the number of DOFs corresponding to a
fine discretization on $\kappa$ is $N_{g}/M$. Hence the total number
of basis functions is $N_{\mc{K}}=n_{b} M$. We also assume $n_b$ is
bounded by a constant while $M,N_{g}$ can increase. In scenario 1, we
increase $N_{g}$ and fix $M$. In scenario 2, $N_{g}$ is proportional
to $M$ while the ratio $N_{g}/M$ is fixed.

The computational cost of the matrix vector multiplication associated
with applying $f(H)$ to $n_{b}$ random vectors is $N_{\text{pole}}
N_{\text{it}} c_{H} n_{b}$. Here $N_{\text{pole}}=r$ is the number of
poles in the rational approximation, $N_{\textit{it}}$ is the number
of iterations to solve for each pole, and $c_{H}$ is the cost of per
iteration.   Since $f(H)$ is a smooth function, $N_{\text{pole}}$
is bounded by a constant independent of $N_{g},M,n_b$. When a good
preconditioner is available, $N_{\textit{it}}$ can also be bounded
by a constant.  $c_{H}$ often is dominated by the matrix-vector
multiplication associated with $H$. Furthermore, when $V$ is a
local potential and when the planewave basis set is used, the cost of
applying $H$ is dominated by applying the Laplacian operator which can
be performed using the fast Fourier transform (FFT). Then $c_{H}\sim
\Or(N_{g}\log N_{g})$.  Since $n_b$ is fixed, in both scenarios
the cost of the matrix-vector multiplication is $\Or(N_{g}\log
N_{g})$.  The cost of each SVD in step~\ref{step:svd_gcalb} is
$\Or((N_{g}/M)\times (n_b+c)^2)$, and the cost for all SVDs is
$\Or(N_{g} (n_b+c)^2)$. Hence the overall complexity for constructing
the GC-ALB set is $\Or(N_{g}\log N_{g})$.  Note that the LC-ALB
approach uses a domain decomposition method, and the computational
complexity of is trivially $\Or(N_{g})$. However, GC-ALB removes
redundant calculations due to overlapping extended elements, and our
numerical results indicate that the efficiency of the GC-ALB approach
can be comparable or even faster when compared to LC-ALB.

The use of the GC-ALB set can also significantly reduce the storage
cost for the spectral projector $P$.  The storage cost for the GC-ALB
set is $(N_{g} / M)\times n_b\times M = N_{g} n_b$.  Viewed as a
matrix, the storage cost for $P$ is $N_{g}^2$. This is generally
very expensive, and $P$ is usually stored using a low rank format
as $P=\Psi\Psi^{*}$, where $\Psi$ is of size $N_{g}\times n$. Then
the storage cost for the coefficient matrix $c_{\kappa,j;i}$ as in
Eq.~\eqref{eqn:DGeig} is $n_{b}M n$, and the total storage cost for
representing $P$ in the GC-ALB set is $N_{g}n_{b} + n_{b} M n$. Hence
when the rank of the projector satisfies
\[
n > \frac{n_{b}}{1-n_{b}M/N_{g}},
\]
the use of the GC-ALB set leads to reduction in the storage cost
for $\Psi$. In practical applications such as Kohn-Sham equations,
this condition is easy to satisfy since $n$ increases with respect
to the system size, while $n_{b}$ is usually a constant on the order
of $10\sim 100$.

Similarly, the computational cost for $\Psi$ using standard
iterative eigensolvers is asymptotically dominated by the need of
orthonormalizing $\Psi$ when $n$ is large. The complexity of this
orthonormalization step scales as $\Or(N_{g} n^2)$. In the GC-ALB
set is constructed, the cost for the orthonormalization is reduced
to $\Or(n_{b} M n^2)$.

Alg.~\ref{alg:gcalb} for constructing the GC-ALB set can also be
efficiently parallelized. For the computation of $W=AR$, which is often
the most time consuming step, the solution of the $N_{\text{pole}}$
shifted linear systems for each column of $R$ are all independent
of each other.  Therefore, the computation can be embarrassingly
parallelized up to $N_{\text{pole}} n_b$ processors. If more than
$N_{\text{pole}}n_b$ processors are available, they will be organized
into $N_{\text{pole}}n_b$ processor groups, and the application of $H$ can
be parallelized as well. For the second part of Alg.~\ref{alg:gcalb},
all calculations can be carried out independently on each element.

\subsection{Generalization to nonlocal potentials}

Another advantage of the GC-ALB approach is that it handles local
and nonlocal potentials on the same footing. The need of computing
spectral projectors associated with nonlocal potentials arise,
for instance, in solving the Hartree-Fock-like equations in quantum
chemistry~\cite{SzaboOstlund1989,Martin2004}. The Hartree-Fock-like
equations require the self-consistent computation of the projector
\begin{equation}
    H[P]=-\frac12 \Delta  +
    V_{\text{ion}} +
    V_{\text{Hxc}}[P] + V_{X}[P], \quad P = \mathbf{1}_{\mc{I}}(H[P]).
  \label{eqn:HF}
\end{equation}
Here the interval $\mc{I}$ contains the lowest $n$ eigenvalues of
$H[P]$. $V_{\text{ion}}, V_{\text{Hxc}}[P]$  are local potentials, and
$V_{X}[P]$ is an integral operator with a nonlocal kernel. Here $[P]$
indicates the nonlinear dependence with respect to $P$. There is no
natural way to consistently incorporate the nonlocal term $V_{X}[P]$
in the LC-ALB approach, while GC-ALB only requires the matrix-vector
multiplication associated with $V_{X}[P]$. A detailed example of
Eq.~\eqref{eqn:HF} will be given in section~\ref{subsec:nonlocal}.

\section{Numerical examples}\label{sec:numer}

We demonstrate the effectiveness of the GC-ALB method for finding
the spectral projector for a linear problem in one, two and three
dimensions in section~\ref{subsec:local}, and for a nonlinear problem
in one dimension in section~\ref{subsec:nonlocal}.  Numerical examples
are performed on Stanford Sherlock cluster bigmem node with quad
socket Intel(R) Xeon(R) CPU E5-4640 @ 2.40GHz and 1.5 TB RAM.  In all
numerical examples, we assume the global domain $\Omega$ satisfies
the periodic boundary condition.  The pseudo-spectral discretization
(a.k.a.  the planewave basis set) provides the reference solution
to the spectral projector, as well as the discretized operator for
performing the matrix-vector multiplication on the global domain in
order to construct the GC-ALB set.  We measure the accuracy of the DG
based methods in terms of the relative error of the eigenvalues within
the range of the spectral projector compared to the reference solution,
defined as
\[
    \frac{\sum_{i\in \mc{I}} \abs{ \varepsilon^{\VN}_{i} -
    \varepsilon_{i}}}{ \sum_{i\in\mc{I}} \abs{ \varepsilon_{i} }}.
\]

The pseudo-spectral discretization can be identified with a set of
uniform grid to discretize $\Omega$.  The integrals needed to construct
the DG bilinear form is done using the Legendre-Gauss-Lobatto (LGL)
grid. A Fourier interpolation procedure is used to interpolate
functions from the uniform grid to the LGL grid, and a stable
barycentric Lagrange interpolation~\cite{BerrutTrefethen2004}
procedure is used to interpolate functions from the LGL grid back
to the uniform grid when needed. In the rational approximation for
the matrix vector multiplication, 16 poles on the upper half complex
plain are actually solved. Since the potential function $V(x)$ in all
numerical examples are real, the rest of the 16 poles are evaluated
via the complex conjugation as in Eq. \eqref{eqn:MZoloFuncMat}. For
each pole, we use the GMRES~\cite{SaadSchultz1986} method to solve the
associated equations with 30 being the restarting number and $10^{-12}$
being the tolerance. The preconditioner is the inverse of a shifted
Laplacian~\cite{GijzenErlanggaVuik2007} with the pole being the shift,
which can be carried out efficiently using fast Fourier transforms
(FFT).  The oversampling parameter $c$ in Alg.~\ref{alg:gcalb} is
set to be 5.  All pseudo-spectral discretized systems, including the
systems for the reference solutions and the system on each extended
element in LC-ALB method, are solved via the LOBPCG~\cite{Knyazev2001}
method, and the associated tolerance is $10^{-12}$ measured in terms
of the maximal residual norm. We use the interior penalty formulation
to patch the discontinuous basis functions to approximate the
eigenfunctions, and the penalty parameter is determined automatically
by solving a local eigenvalue problem as in~\cite{LinStamm2016}.

\subsection{Linear problems with local potentials}\label{subsec:local}

\subsubsection{One dimensional case}

\begin{figure}[h]
  \begin{center}
    \subfloat[]{\includegraphics[width=0.44\textwidth]{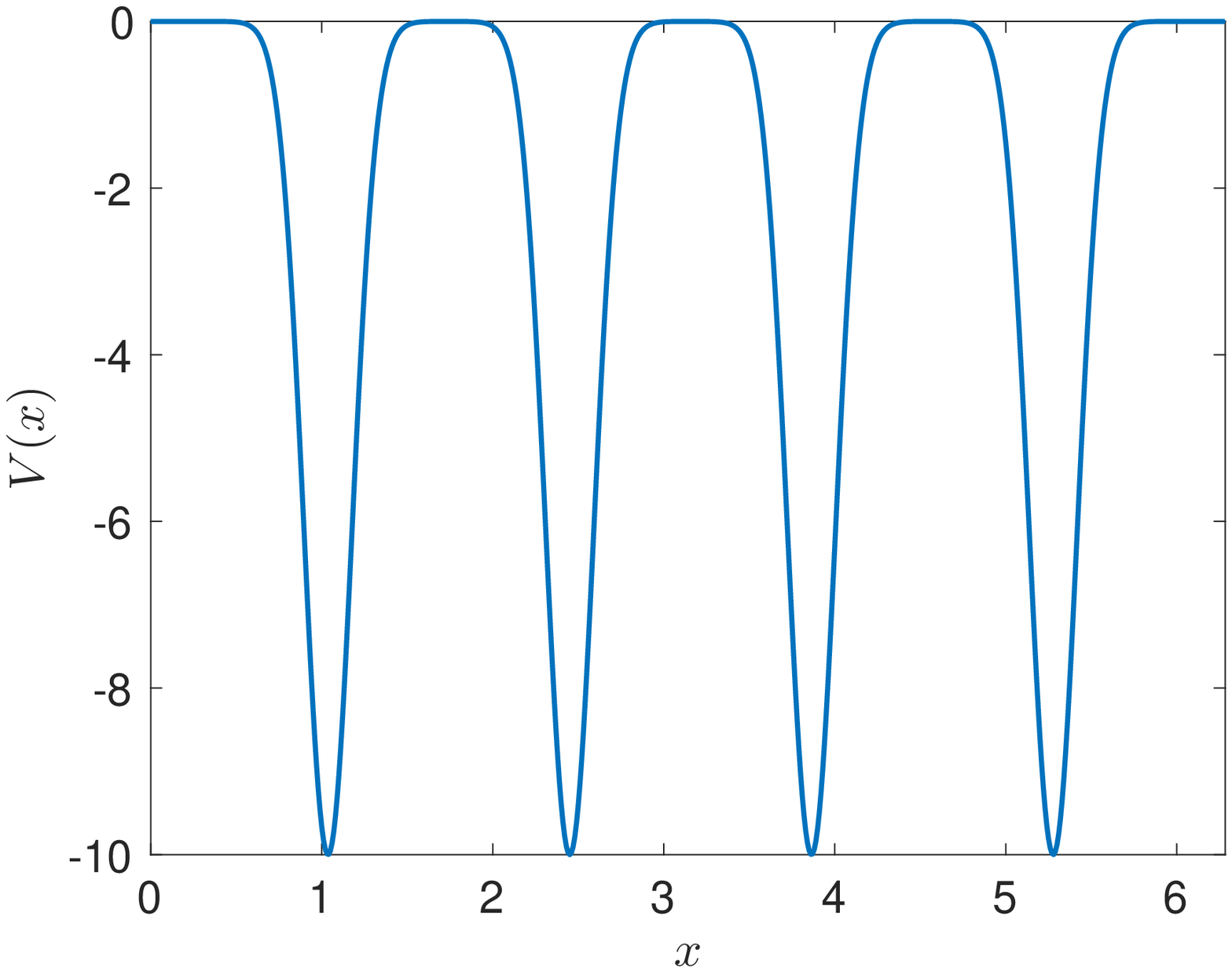}}~~
    \subfloat[]{\includegraphics[width=0.45\textwidth]{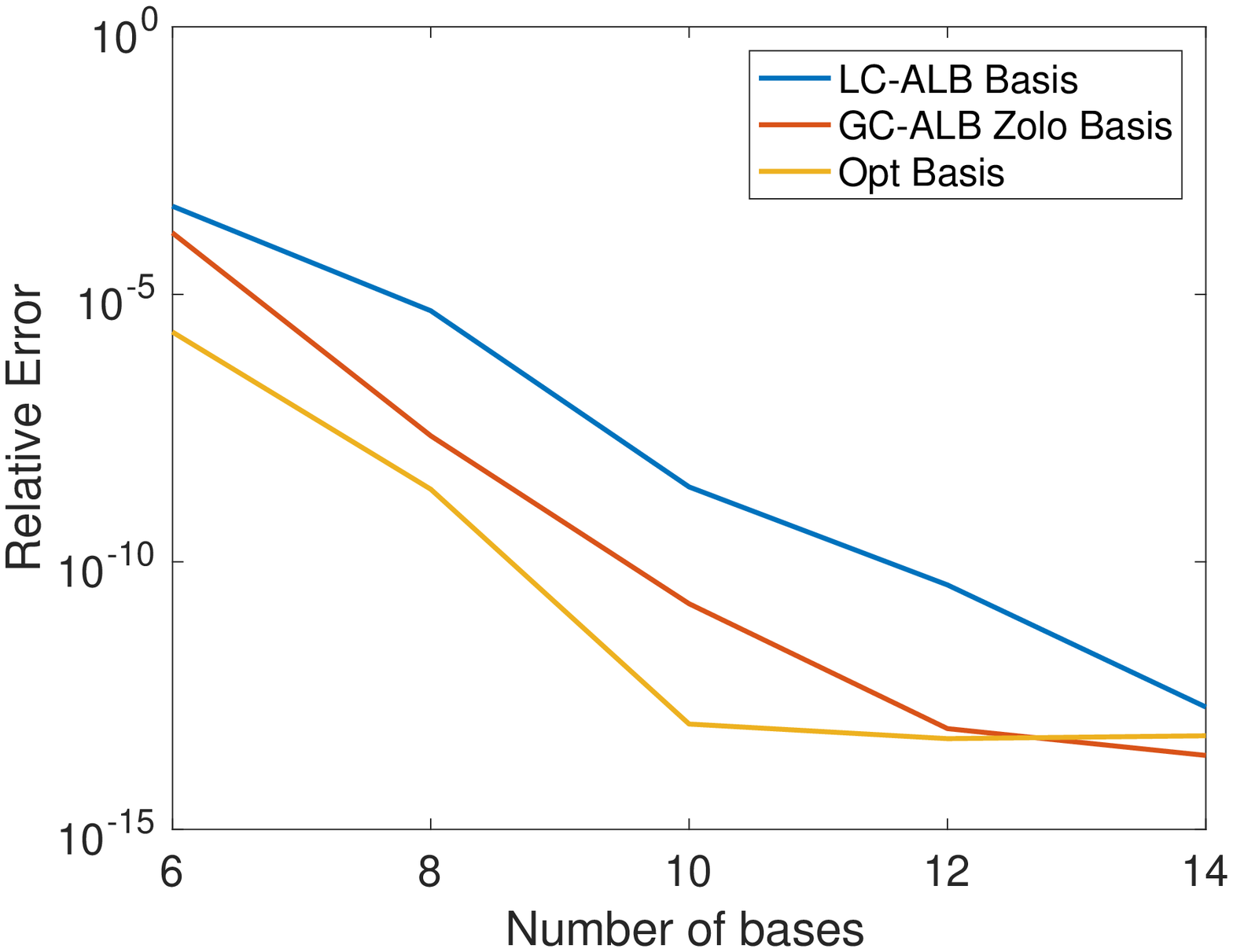}}
  \end{center}
    \caption{(a) the 1D potential function (b) the relative
    errors of LC-ALB, GC-ALB, and Opt method for different number
    of basis functions.}
  \label{fig:NumResLoc1D}
\end{figure}

Our first example is a second order differential
operator~\eqref{eqn:Hoperator} on $\Omega = (0,2\pi)$ in 1D.  $V$
is a local potential with four Gaussian potential wells at positions
$x = \{ 1.0367, 2.4504, 3.8642, 5.2779 \}$. The depth for each
well is $-10.0$ whereas the standard deviation is set to be $0.2$.
Fig.~\ref{fig:NumResLoc1D}~(a) shows the potential $V(x)$. The interval
$\mc{I}$ associated with the spectral projector $P$ is assumed to
cover the lowest 16 eigenvalues.

The global domain $\Omega$ is partitioned into $7$ elements. Within
each element, $40$ LGL grid points are used to evaluate the
integrals in the DG bilinear form accurately.  The pseudo-spectral
method discretizes $\Omega$ using $140$ planewave basis functions,
which can be identified with a uniform grid with $140$ grid points.
Under these settings, three adaptive local basis construction methods
are considered, i.e., LC-ALB, GC-ALB with rational approximation
for the projector (GC-ALB), and the optimal basis set (Opt).
For different methods, we vary the number of basis functions used
in each element from 6 to 14. The relative error of the smallest
16 eigenvalues is measured against a reference solution, which is
calculated via the pseudo-spectral method with 500 planewave basis
functions directly.  In the GC-ALB, we set the interval as $\mc{I} =
[a,b] = [\lambda_1,\lambda_{16}]$ and the gap parameters as $a_-
= -\infty$ and $b_+ = \lambda_{16}+1.0$, where $\lambda_1$ and
$\lambda_{16}$ denote the smallest eigenvalue and the 16th smallest
eigenvalue, respectively.  Fig.~\ref{fig:NumResLoc1D}~(b) shows the
relative errors for different methods with varying number of basis
functions. More details are reported in Tab.~\ref{tab:NumResLoc1D}.

\begin{table}[h]
    \centering
    \begin{tabular}{cccccc}
\toprule
 Method & $n_b$ & err & $T_{Basis}$ (sec) & $T_{DG}$ (sec) & $n_{\text{tot iter}}$ \\
\toprule
\multirow{5}{*}{\shortstack{GC-ALB}} & 6 & 1.41e-04 & 5.69e-01 & 1.12e-02 & 147 \\
 & 8 & 2.27e-08 & 6.14e-01 & 1.06e-02 & 147 \\
 & 10 & 1.65e-11 & 6.22e-01 & 1.18e-02 & 148 \\
 & 12 & 7.64e-14 & 7.21e-01 & 1.49e-02 & 149 \\
 & 14 & 2.41e-14 & 7.24e-01 & 1.26e-02 & 150 \\
\toprule
\multirow{5}{*}{LC-ALB} & 6 & 4.45e-04 & 2.91e-01 & 1.05e-02 & - \\
 & 8 & 4.93e-06 & 3.24e-01 & 1.06e-02 & - \\
 & 10 & 2.51e-09 & 3.59e-01 & 1.55e-02 & - \\
 & 12 & 3.71e-11 & 3.91e-01 & 1.42e-02 & - \\
 & 14 & 1.93e-13 & 4.20e-01 & 1.25e-02 & - \\
\bottomrule
    \end{tabular}
  \caption{Numerical results for GC-ALB method and LC-ALB method. $n_b$
    is the maximum number of basis functions for each element, err is the
    relative error of the smallest 16 eigenvalues, $T_{Basis}$ and $T_{DG}$
    are the runtime for basis construction and DG solving respectively,
    $n_{\text{tot iter}}$ is the total number of iterations for solving
    linear systems throughout the algorithm.}
  \label{tab:NumResLoc1D}
\end{table}

For the one dimensional operator, as shown in
Fig.~\ref{fig:NumResLoc1D}~(b), the relative errors for all three
methods decay exponentially as the number of basis functions increases.
As discussed in section~\ref{sec:Optbasis}, the Opt basis defines
the optimal discontinuous basis set for a given partition of the
global domain and number of basis functions in each element, and this
is confirmed in Fig.~\ref{fig:NumResLoc1D}~(b). On the other hand,
the performance both GC-ALB and LC-ALB closely follow the Opt basis.
Given the same number of basis functions, GC-ALB is about one digit
more accurate than LC-ALB. When the number of basis functions is larger
than or equal to 14, both methods reach the numerical accuracy limit
and can not be further improved. The runtime of the GC-ALB method
and LC-ALB are about the same. The numbers of total iterations are
about 148, which means the iteration number for solving each pole in
\eqref{eqn:MZoloFuncMat} is on average smaller than 10.

\subsubsection{Two dimensional case}

\begin{figure}[h]
  \begin{center}
    \subfloat[]{\includegraphics[width=0.44\textwidth]{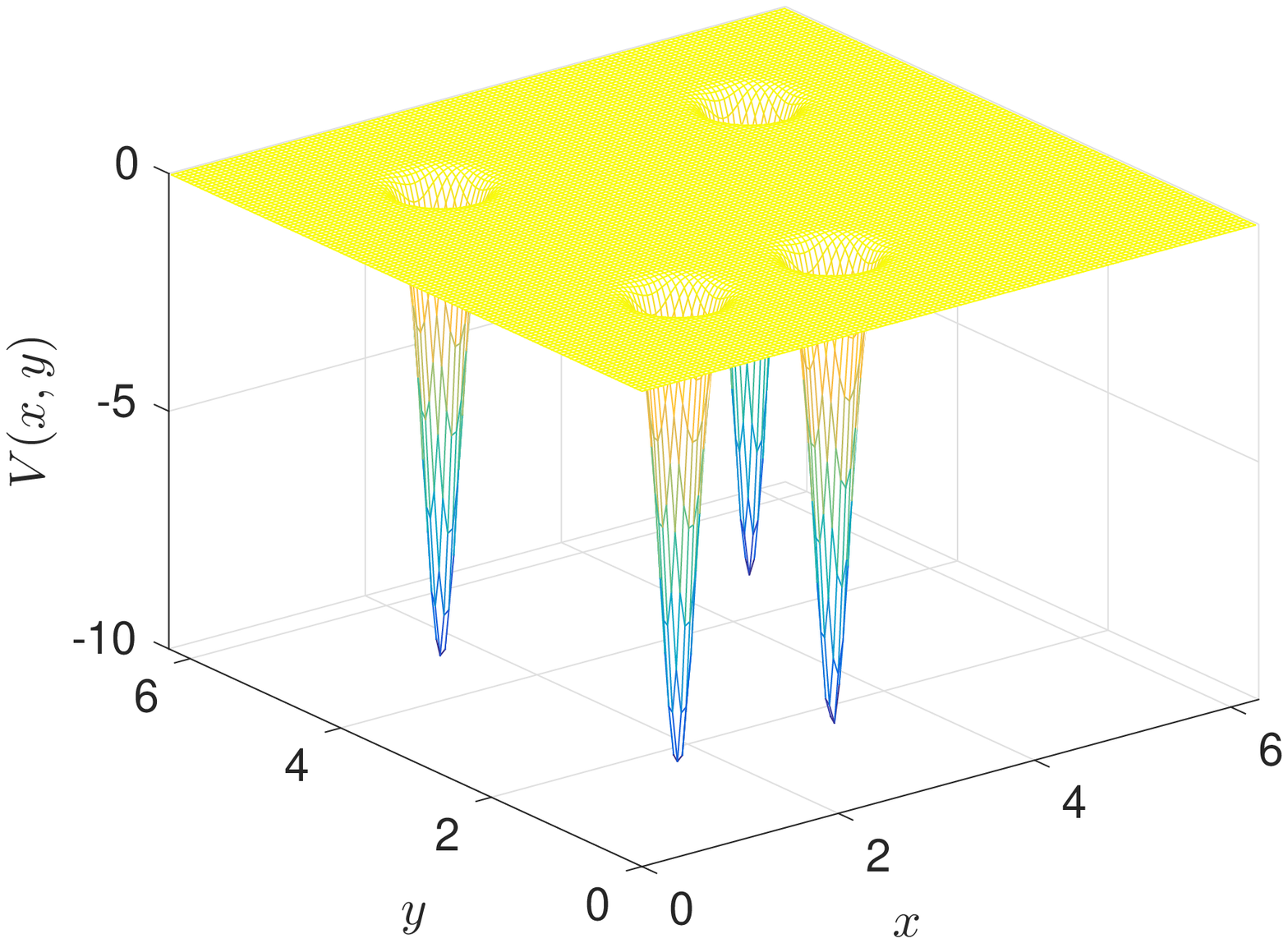}}~~
    \subfloat[]{\includegraphics[width=0.45\textwidth]{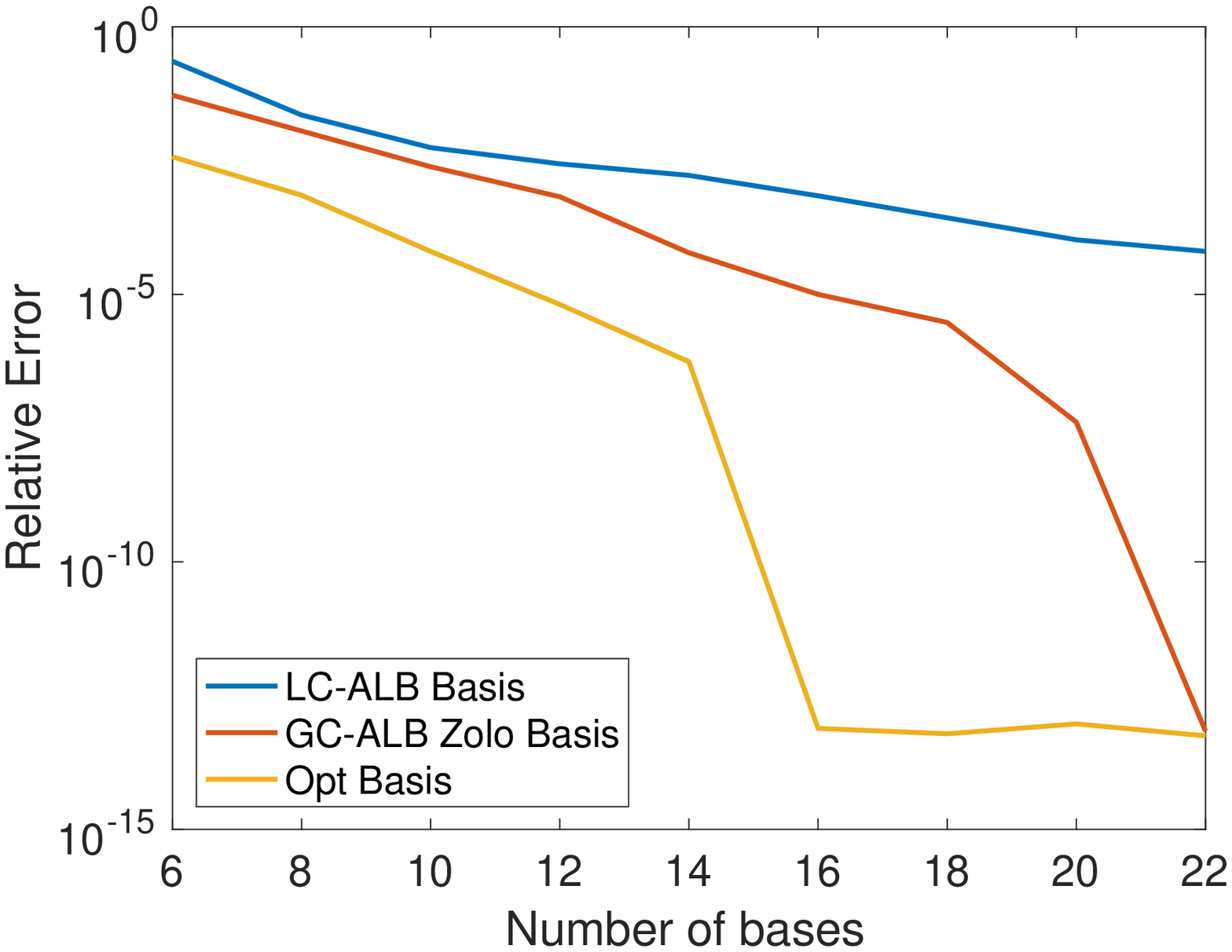}}
  \end{center}
    \caption{(a) the 2D potential function (b) the relative
    errors of LC-ALB, GC-ALB, and Opt method for different
    number of basis functions.}
  \label{fig:NumResLoc2D}
\end{figure}

This example is a second order differential
operator~\eqref{eqn:Hoperator} on $\Omega = (0,2\pi)^2$ in 2D.
$V$ is a local potential with four Gaussian wells as shown in
Fig.~\ref{fig:NumResLoc2D}~(a). The depth for each well is $-10.0$ and
the standard deviation is $0.2$. Similar to one dimensional example,
the interval $\mc{I}$ associated with the spectral projector $P$
is assumed to cover the lowest 16 eigenvalues.

The global domain $\Omega$ is partitioned into $7 \times 7$
elements. Within each element, $40 \times 40$ two dimensional LGL
grid points are used to evaluate the integrals in the DG bilinear
form accurately.  The pseudo-spectral method discretizes $\Omega$
using $140^2$ planewave basis functions, which can be identified with
a uniform two dimensional grid with $140 \times 140$ grid points.
Similar name conventions for LC-ALB, GC-ALB and Opt are used as in the
one dimensional example.  For different methods, we vary the number
of basis functions used in each element from 8 to 22. The relative
error of the smallest 16 eigenvalues is measured against a reference
solution, which is calculated via the pseudo-spectral method with
$300^2$ planewave basis functions directly.  In the GC-ALB, we set
the interval as $\mc{I} = [a,b] = [\lambda_1,\lambda_{16}]$ and the
gap parameters as $a_- = -\infty$ and $b_+ = \lambda_{16}+0.1$, where
$\lambda_1$ and $\lambda_{16}$ denote the smallest eigenvalue and the
16th smallest eigenvalue.  Fig.~\ref{fig:NumResLoc2D}~(b) shows the
relative errors for different methods with varying number of basis
functions. More details are reported in Tab.~\ref{tab:NumResLoc2D}.

\begin{table}[h]
    \centering
    \begin{tabular}{cccccc}
\toprule
 Method & $n_b$ & err & $T_{Basis}$ (sec) & $T_{DG}$ (sec) & $n_{\text{tot iter}}$ \\
\toprule
\multirow{9}{*}{\shortstack{GC-ALB}} & 6 & 5.27e-02 & 2.72e+01 & 2.53e-01 & 156 \\
 & 8 & 1.13e-02 & 2.95e+01 & 2.80e-01 & 156 \\
 & 10 & 2.43e-03 & 3.37e+01 & 3.29e-01 & 156 \\
 & 12 & 6.67e-04 & 3.71e+01 & 3.95e-01 & 156 \\
 & 14 & 5.98e-05 & 3.83e+01 & 4.76e-01 & 156 \\
 & 16 & 9.99e-06 & 4.15e+01 & 4.61e-01 & 156 \\
 & 18 & 2.98e-06 & 4.12e+01 & 6.31e-01 & 156 \\
 & 20 & 4.10e-08 & 4.59e+01 & 6.51e-01 & 157 \\
 & 22 & 6.87e-14 & 4.99e+01 & 6.78e-01 & 156 \\
\toprule
\multirow{9}{*}{LC-ALB} & 6 & 2.27e-01 & 1.41e+01 & 2.21e-01 & - \\
 & 8 & 2.25e-02 & 1.49e+01 & 2.45e-01 & - \\
 & 10 & 5.53e-03 & 2.59e+01 & 3.12e-01 & - \\
 & 12 & 2.75e-03 & 2.52e+01 & 3.75e-01 & - \\
 & 14 & 1.67e-03 & 3.30e+01 & 4.75e-01 & - \\
 & 16 & 6.95e-04 & 2.76e+01 & 5.17e-01 & - \\
 & 18 & 2.69e-04 & 2.85e+01 & 5.90e-01 & - \\
 & 20 & 1.05e-04 & 2.90e+01 & 6.21e-01 & - \\
 & 22 & 6.37e-05 & 4.38e+01 & 7.39e-01 & - \\
\bottomrule
    \end{tabular}
  \caption{Numerical results of GC-ALB method and LC-ALB method for the
  two dimensional example.}
  \label{tab:NumResLoc2D}
\end{table}

Fig.~\ref{fig:NumResLoc2D}~(b) shows that the differences among LC-ALB,
GC-ALB and Opt basis sets become more significant in 2D.  The relative
errors for Opt and GC-ALB decreases to the level of $10^{-14}$ when 16
and 22 basis functions are constructed for each element respectively.
On the other hand side, the relative errors for LC-ALB method remains
around $6.37\times 10^{-5}$ when $22$ basis functions are used for
each element. In order to achieve an relative error that is below
$10^{-12}$, we also find that $120$ basis functions per element are
needed in the LC-ALB approach. Tab.~\ref{tab:NumResLoc2D} shows that
the cost for the GC-ALB and LC-ALB approaches are comparable in 2D. The
fluctuation of the runtime in the LC-ALB approach is mostly due to
the fluctuation of the number of iterations for the LOBPCG solver.
In the GC-ALB approach, the number of iterations for solving each pole
here is around 9 on average for all cases, which gives $n_{\text{tot
iter}}$ to be around $156$ in all cases.

Below we demonstrate the weak scaling performance of the GC-ALB set
in 2D.  Starting from the potential in Fig.~\ref{fig:NumResLoc2D}
(a), we increase the size of the domain by periodically repeating the
potential along $x$ and $y$ directions by a factor of $\ell$.  We vary
$\ell$ from 1 to 6, as shown in Tab.~\ref{tab:NumResLoc2DWeak}, and
the domain $\Omega$ is extended from $(0,2\pi)^2$ to $(0,12\pi)^2$. The
number of planewave basis functions, the number of elements, as well as
the number of eigenvalues to be computed are proportional to the size
of $\Omega$.  The parameters used within each element are the same as
before and 20 basis functions are constructed for each element. In
terms of the parameters in Zolotarev's function approximation, we
set the interval as $\mc{I} = [a,b] = [\lambda_1,\lambda_{n}]$ and
the gap parameters as $a_- = -\infty$ and $b_+ = \lambda_{n}+0.1$,
where $\lambda_1$ and $\lambda_{n}$ denote the smallest eigenvalue
and the $n$th smallest eigenvalue. The relative error of the smallest
$n$ eigenvalues is measured against reference solutions, which are
calculated via the pseudo-spectral method with $(300\ell)^2$ planewave
basis functions directly. Since the reference solution for $\ell=6$
cannot be finished within the limited runtime on the Sherlock system,
only the GC-ALB runtime and the total iteration number are reported
here.

\begin{table}[h]
    \centering
    \begin{tabular}{ccccccc}
\toprule
 $\ell$ & $\Omega$ & $n$ & err & $T_{Basis}$ (sec) & $T_{DG}$ (sec) & $n_{\text{tot iter}}$ \\
\toprule
 1 & $(0,2\pi)^2$ & 16 & 1.14e-07 & 4.97e+01 & 8.74e-01 & 156 \\
 2 & $(0,4\pi)^2$ & 64 & 4.47e-07 & 3.17e+02 & 6.29e+00 & 197 \\
 3 & $(0,6\pi)^2$ & 144 & 5.78e-07 & 6.21e+02 & 4.64e+01 & 209 \\
 4 & $(0,8\pi)^2$ & 256 & 7.29e-07 & 2.09e+03 & 1.87e+02 & 271 \\
 5 & $(0,10\pi)^2$ & 400 & 6.33e-07 & 3.53e+03 & 6.11e+02 & 268 \\
 6 & $(0,12\pi)^2$ & 576 & - & 7.52e+03 & 1.66e+03 & 268 \\
\bottomrule
    \end{tabular}
  \caption{Numerical results of the weak scaling of the GC-ALB method for
  the two dimensional example. Here 20 basis functions are used within each
  element, $\ell$ denotes the number of repeated domain on each dimension,
  $n$ denotes the number of calculated eigenvalues.}
  \label{tab:NumResLoc2DWeak}
\end{table}

\begin{figure}[h]
  \begin{center}
    \subfloat[]{\includegraphics[width=0.48\textwidth]{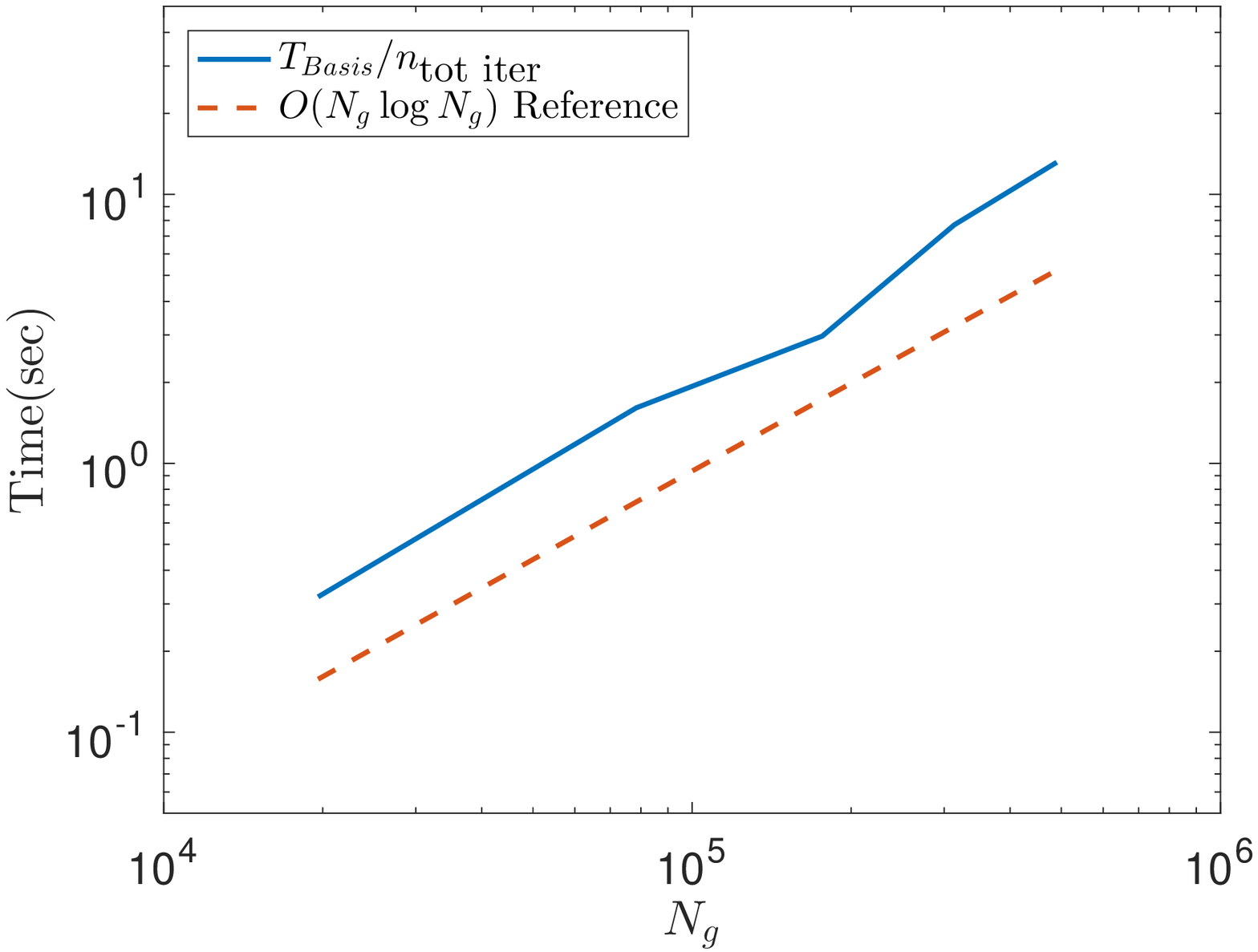}}~~
    \subfloat[]{\includegraphics[width=0.48\textwidth]{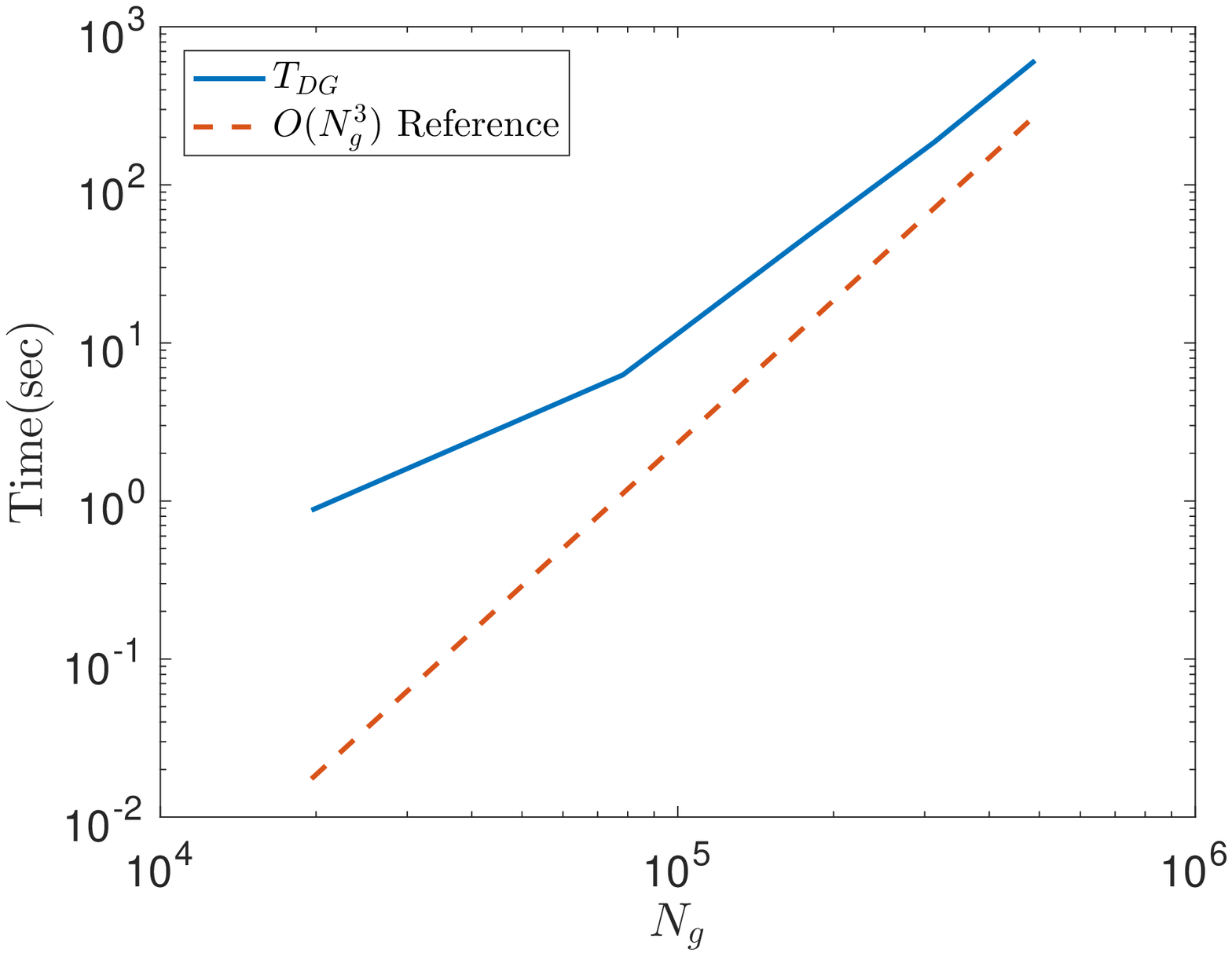}}
  \end{center}
    \caption{Scalings of (a) average single iteration runtime and (b)
    DG solving time.}
  \label{fig:WeakScaling}
\end{figure}

Tab.~\ref{tab:NumResLoc2DWeak} shows that the relative errors are
approximately the same for all $\ell$. The runtime, $T_{Basis}$,
for the basis construction in GC-ALB increases proportional to
$\ell^2\log \ell$, which means $T_{Basis}$ is quasi-linear in the
number of planewave basis functions (see Fig.~\ref{fig:WeakScaling}
(a)). Meanwhile $T_{DG}$, which is the cost for solving the DG
problem is super-linear with respect to the number of planewave basis
functions (see Fig.~\ref{fig:WeakScaling} (b)).  $T_{Basis}$ is consist
with the complexity analysis and $T_{DG}$ is close aligned with the
complexity analysis when $N_g$ is large.  We observe that when $\ell$
is relatively small, the number of total iterations mildly increases
with respect to $\ell$. As $\ell$ keeps on increasing, the total
iteration number stays around $270$.

\subsubsection{Three dimensional case}

This example is a second order differential
operator~\eqref{eqn:Hoperator} on $\Omega = (0,2\pi)^3$ in 3D.  $V$
is a local potential with four Gaussian wells.  The depth for each
well is $-10$ whereas the standard deviation is set to be $0.2$.
Fig.~\ref{fig:NumResLoc2D}~(a) shows the isosurface for the potential
function $V(x,y,z) = -1$. Similar to previous examples, the interval
$\mc{I}$ associated with the spectral projector $P$ is assumed to
cover the lowest 16 eigenvalues.

\begin{figure}[h]
  \begin{center}
    \subfloat[]{\includegraphics[width=0.44\textwidth]{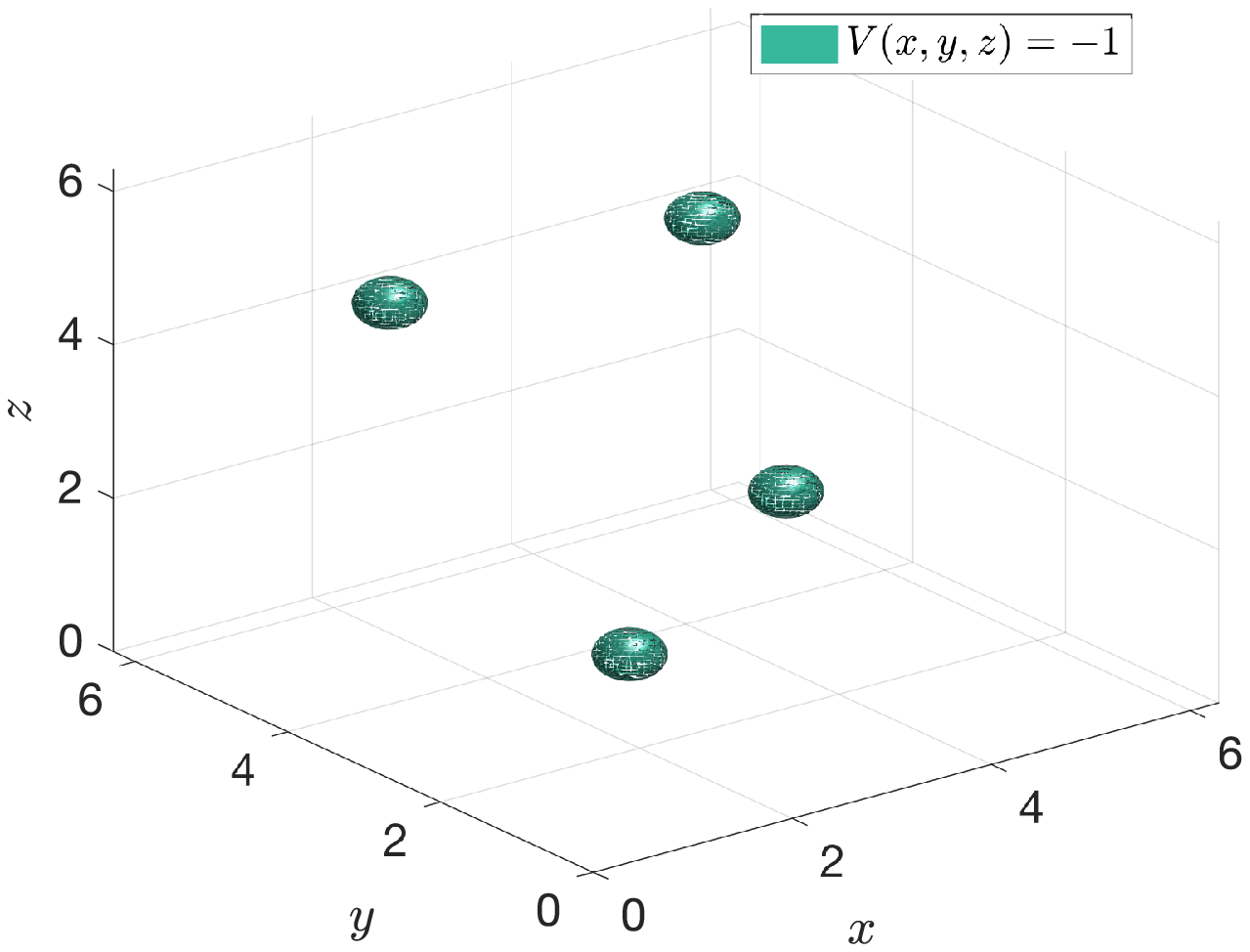}}~~
    \subfloat[]{\includegraphics[width=0.45\textwidth]{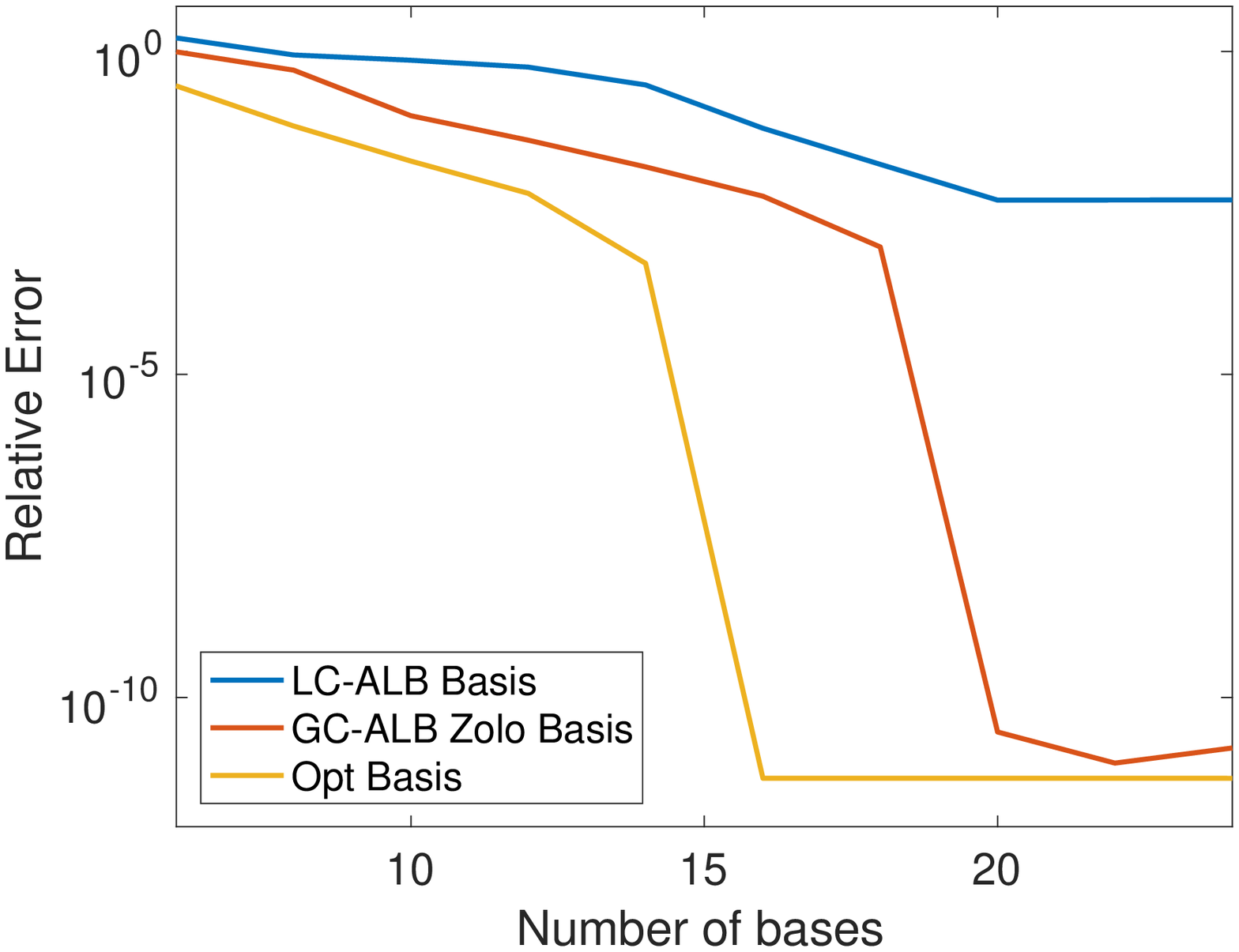}}
  \end{center}
    \caption{(a) isosurface plot $V(x,y,z) = -1$ for the 3D potential
    function (b) the relative error of LC-ALB, GC-ALB and
    Opt method with different numbers of basis functions.}
  \label{fig:NumResLoc3D}
\end{figure}

The global domain $\Omega$ is partitioned into $4 \times 4 \times
4$ elements. Within each element, $30 \times 30 \times 30$ three
dimensional LGL grid points are used to evaluate the integrals in the
DG bilinear form accurately.  The pseudo-spectral method discretizes
$\Omega$ using $60^3$ planewave basis functions, which can be
identified with a uniform three dimensional grid with $60 \times 60
\times 60$ grid points.  Similar name conventions for LC-ALB, GC-ALB
and Opt are used as in the one dimensional example.  For different
methods, we vary the number of basis functions used in each element
from 8 to 24. The relative error of the smallest 16 eigenvalues
is measured against a reference solution, which is calculated via
the pseudo-spectral method with $100^3$ planewave basis functions
directly.  In the GC-ALB, we set the interval as $\mc{I} = [a,b] =
[\lambda_1,\lambda_{16}]$ and the gap parameters as $a_- = -\infty$
and $b_+ = \lambda_{16}+0.01$, where $\lambda_1$ and $\lambda_{16}$
denote the smallest eigenvalue and the 16th smallest eigenvalue.
Fig.~\ref{fig:NumResLoc3D}~(b) shows the relative errors for different
methods with varying number of basis functions. More details are
reported in Tab.~\ref{tab:NumResLoc3D}.

\begin{table}[h]
    \centering
    \begin{tabular}{cccccc}
\toprule
 Method & $n_b$ & err & $T_{Basis}$ (sec) & $T_{DG}$ (sec) & $n_{\text{tot iter}}$ \\
\toprule
\multirow{10}{*}{\shortstack{GC-ALB}} & 6 & 9.89e-01 & 1.52e+02 & 3.07e+00 & 128 \\
 & 8 & 5.15e-01 & 1.60e+02 & 3.60e+00 & 128 \\
 & 10 & 1.01e-01 & 2.07e+02 & 4.34e+00 & 129 \\
 & 12 & 4.24e-02 & 1.81e+02 & 4.87e+00 & 129 \\
 & 14 & 1.64e-02 & 2.43e+02 & 5.68e+00 & 129 \\
 & 16 & 5.76e-03 & 2.40e+02 & 6.48e+00 & 129 \\
 & 18 & 9.42e-04 & 2.89e+02 & 7.66e+00 & 129 \\
 & 20 & 2.92e-11 & 3.04e+02 & 8.30e+00 & 129 \\
 & 22 & 9.63e-12 & 3.19e+02 & 8.84e+00 & 129 \\
 & 24 & 1.65e-11 & 3.51e+02 & 1.23e+01 & 129 \\
\toprule
\multirow{10}{*}{LC-ALB} & 6 & 1.62e+00 & 5.25e+02 & 3.15e+00 & - \\
 & 8 & 8.84e-01 & 2.02e+03 & 3.31e+00 & - \\
 & 10 & 7.32e-01 & 1.54e+03 & 5.23e+00 & - \\
 & 12 & 5.73e-01 & 1.22e+03 & 6.18e+00 & - \\
 & 14 & 3.04e-01 & 1.16e+03 & 7.05e+00 & - \\
 & 16 & 6.48e-02 & 1.26e+03 & 7.82e+00 & - \\
 & 18 & 1.79e-02 & 1.46e+03 & 9.00e+00 & - \\
 & 20 & 5.00e-03 & 1.87e+03 & 9.19e+00 & - \\
 & 22 & 5.02e-03 & 1.81e+03 & 1.07e+01 & - \\
 & 24 & 5.04e-03 & 1.96e+03 & 1.20e+01 & - \\
\bottomrule
    \end{tabular}
  \caption{Numerical results of GC-ALB method and LC-ALB method for the
  three dimensional example.}
  \label{tab:NumResLoc3D}
\end{table}

\begin{table}[h]
    \centering
    \begin{tabular}{cccc|cccc}
\toprule
        & $n_b$ & err & DOFs & & $n$ & err & DOFs \\
\toprule
 \multirow{4}{*}{\shortstack{GC-ALB}} & 14 & 1.69e-02 &  896 &
 \multirow{4}{*}{planewave} & 16 & 9.50e-03 & 4096  \\
 & 16 & 4.34e-03 & 1024 & & 20 & 2.13e-03 & 8000 \\
 & 18 & 1.36e-04 & 1152 & & 26 & 1.60e-04 & 17576 \\
 & 20 & 9.19e-12 & 1280 & & 60 & 1.10e-10 & 216000 \\
\bottomrule
    \end{tabular}
  \caption{Comparison of the degrees of freedom for GC-ALB method and
    planewave method for the three dimensional example.}
  \label{tab:NumResLoc3DComp}
\end{table}

For three dimensional systems, the GC-ALB method exhibits even clearer
advantage over the LC-ALB method.  In Fig.~\ref{fig:NumResLoc3D}~(b)
and Tab~\ref{tab:NumResLoc3D}, the relative errors for GC-ALB method
decay quickly to the level of $10^{-11}$, while the asymptotic decay
rate of the LC-ALB method is much slower. The GC-ALB approach is
also more efficient in terms of the runtime. For most of the cases in
Tab~\ref{tab:NumResLoc3D}, the GC-ALB method is about 6 times faster
than LC-ALB method. The numbers of the applications of the operator
to test vectors are 129 in GC-ALB method for all different number
of bases.  In addition, Tab.~\ref{tab:NumResLoc3DComp} shows that the
number of DOFs for the GC-ALB set is much smaller than that needed
for the planewave basis set to reach the same level of accuracy.
Here the DOFs for the GC-ALB set is equal to the dimension of the
DG matrix, and the DOFs for the planewave basis set is the number of
planewave basis functions.

\subsection{Nonlinear problems with nonlocal
potentials}\label{subsec:nonlocal}

In order to test the effectiveness of the GC-ALB approach for nonlocal
potentials, we consider the following model for Hartree-Fock-like
equations in one dimension. The Hamiltonian operator acting on a
function $\psi$ is given by
\begin{equation}
  \begin{split}
   (H[P]\psi)(x) = &-\frac{1}{2} \frac{d^2}{dx^2}\psi(x) + \left(\int
   K(x,y) (m(y)+P(y,y)) \ud y \right)\psi(x) \\ &- \alpha \int K(x,y)
   P(x,y) \psi(y) \ud y
  \end{split}
  \label{eqn:HrHF}
\end{equation}
Compared to Eq.~\eqref{eqn:HF}, the second term on the right
hand side of Eq.~\eqref{eqn:HrHF} corresponds to $V_{\text{ion}}$
and $V_{\text{Hxc}}[P]$ and is a local potential, while the third
term corresponds to $V_{X}[P]$ and is a nonlocal potential.  Here
$m(x)=\sum_{i=1}^{M} m_i(x-R_i)$, with the position of the $i$-th
nuclei denoted by $R_{i}$.  Each function $m_{i}(x)$ takes the form
\begin{equation} 
  m_{i}(x) = -\frac{Z_i}{\sqrt{2\pi\sigma_{i}^2}}
  e^{-\frac{x^2}{2\sigma_i^2}},
  \label{}
\end{equation}
where $Z_i$ is an integer representing the charge of the $i$-th
nucleus.  Instead of using a bare Coulomb interaction, which diverges
in 1D, we adopt a Yukawa kernel
\begin{equation}
  K(x,y) = \frac{2\pi e^{-\mu \abs{x-y}}}{\mu\epsilon_{0}},
  \label{eqn:YukawaK1}
\end{equation}
which satisfies the equation
\begin{equation}
  -\frac{d^2}{d x^2} K(x,y) + \mu^2 K(x,y) = \frac{4\pi}{\epsilon_0}
  \delta(x-y).  \label{eqn:YukawaK2}
\end{equation}
As $\mu\to 0$, the Yukawa kernel approaches the bare Coulomb
interaction given by the Poisson equation. The parameters
$\epsilon_0,\alpha$ are used to ensure that the contribution from
different terms are comparable. And the notations here are different
from the ones in Section~\ref{sec:RationalApprox}. In this example,
we choose $\Omega=(0,80)$, $M=8$, $\sigma_{i}=3.0$, $Z_{i}=2.0$,
$\mu=0.01$, $\epsilon_{0}=10$, $\alpha=0.05$.  Besides these
parameters, for the Zolotarev's function approximation in every
iteration, $16$ poles are used, $a_-=-\infty$, $a$ is the smallest
eigenvalue calculated each iteration, $b = -3.388$ which is the
converged Fermi level, and $b_+ = 0$.  The self-consistent spectral
projector $P$ is given by the lowest $16$ eigenfunctions of $H[P]$.

In order to find the self-consistent spectral projector,
we use a two level self consistent field (SCF) iteration
that is commonly adopted to solve such Hartree-Fock-like
equations~\cite{GiannozziBaroniBoniniEtAl2009,Lin2016ACE}.  The SCF
iterations are split into an outer loop and an inner loop.  At the
beginning of each outer SCF loop, we update the nonlocal potential
$V_{X}[P]$ using a fixed point iteration, i.e. $P$ is updated by the
converged spectral projector $P$ from the inner SCF loop.  In the
inner SCF loop, we fix the nonlocal potential $V_{X}[P]$ as if it
were independent of $P$, and update the local potential via the
diagonal part of the projector $P(x,x)$ using the Anderson mixing
method for charge mixing~\cite{Anderson1965}. The convergence of the
outer iteration is measured by the convergence of the exchange energy
defined as
\begin{equation}
  E_{X} = -\int P(x,y)K(x,y)P(x,y) \ud x \ud y.
  \label{eqn:exchangeenergy}
\end{equation}

In each inner SCF iteration, we apply the GC-ALB method with
Zolotarev's function approximation together with DG method to construct
the spectral projector efficiently, which is denoted as ``GC-ALB''
in the rest of this paper. As a comparison, we also conduct the inner
and outer SCF iterations with the spectral projector calculated via
planewave method, which is denoted as ``planewave''.

\begin{figure}[h]
  \begin{center}
    \subfloat[]{\includegraphics[width=0.4\textwidth]{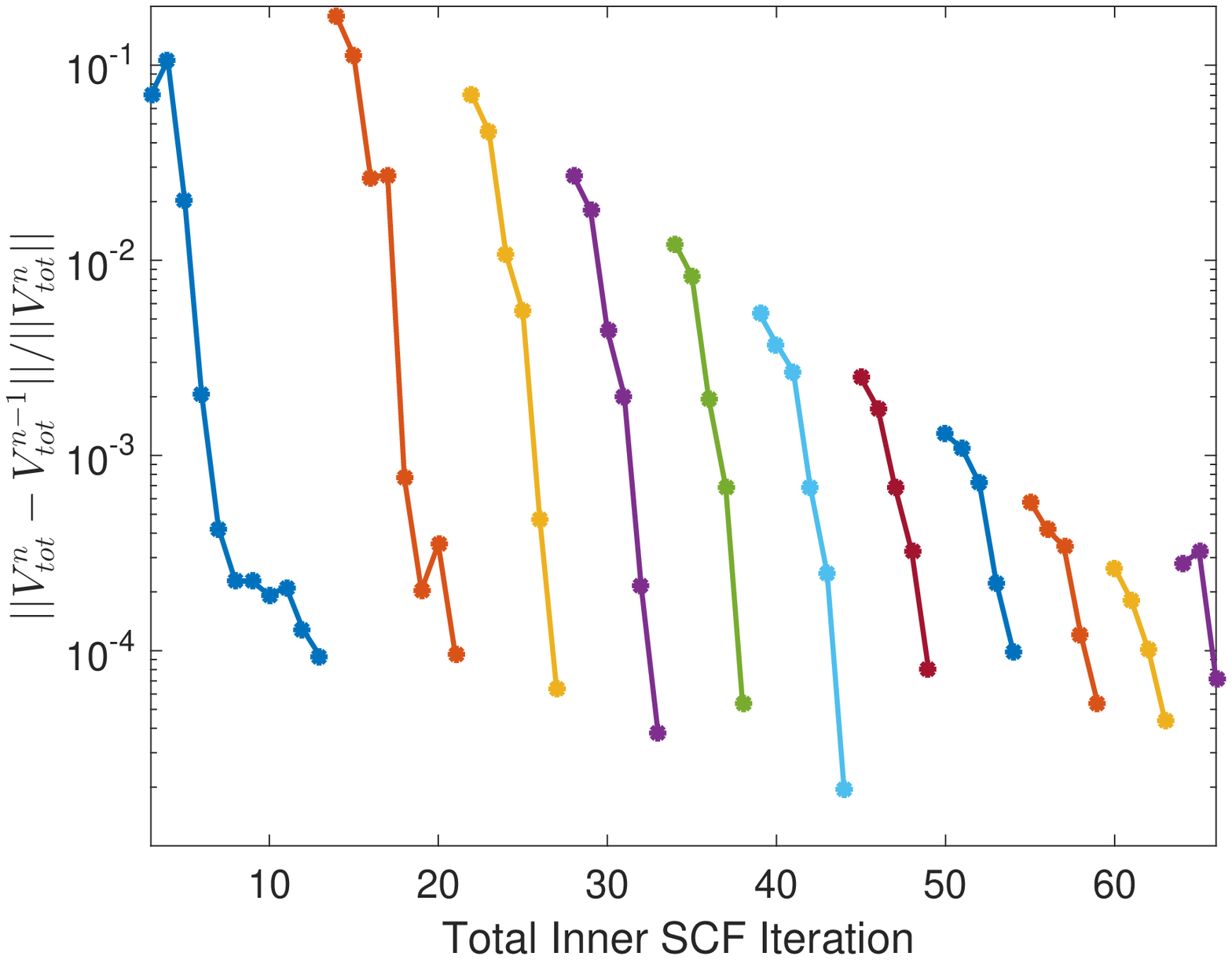}}~~
    \subfloat[]{\includegraphics[width=0.4\textwidth]{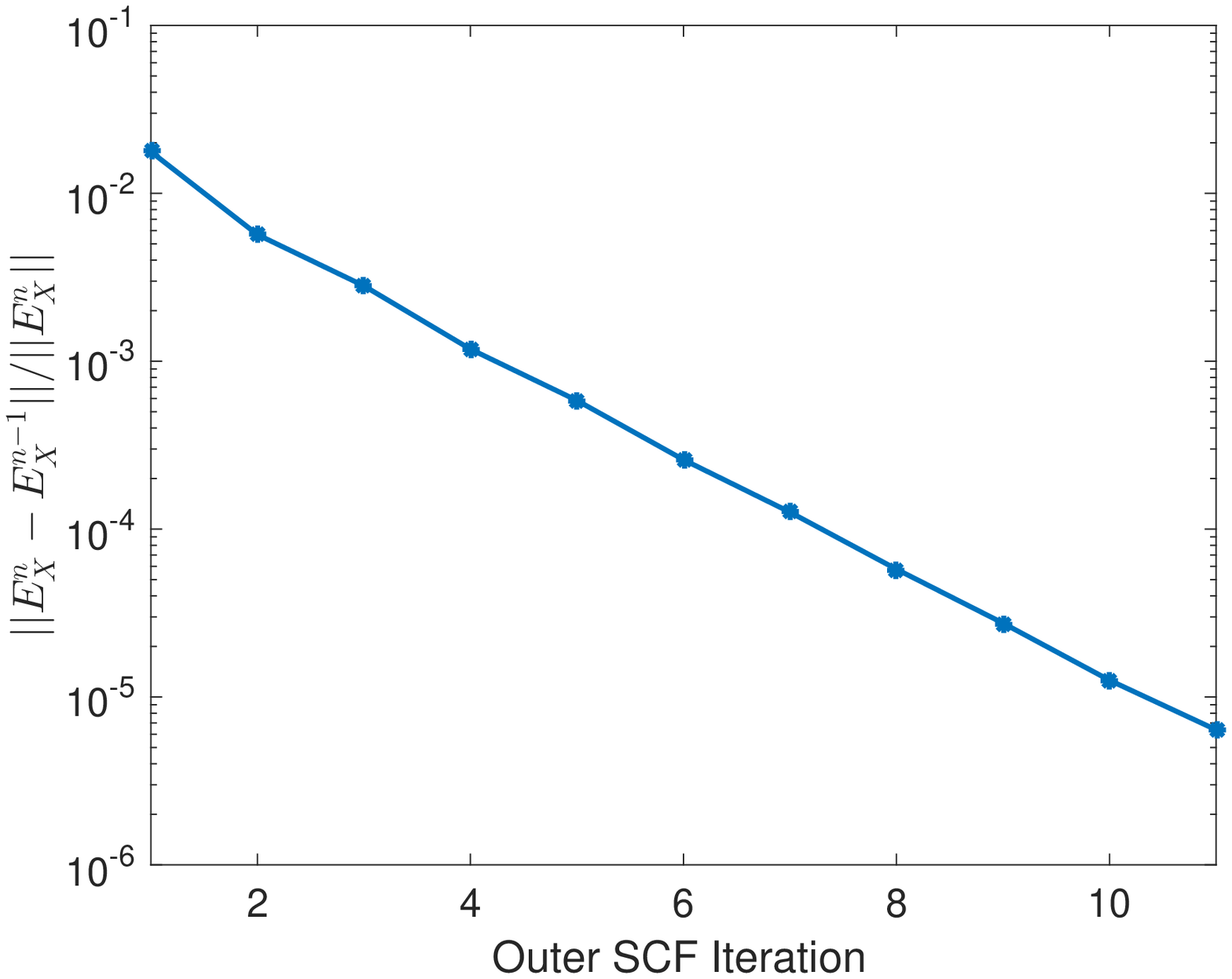}}
  \end{center}
  \caption{GC-ALB: (a) The relative errors of the total local potential.
    Each point is a inner SCF iteration whereas each color line indicates a
    outer SCF iteration. (b) The relative errors of the energy associated
    with the nonlocal potential.} \label{fig:NumResSCF}
\end{figure}

Fig.~\ref{fig:NumResSCF} (a) and (b) show the convergence
behavior of the two level SCF iterations using the GC-ALB
set. Fig.~\ref{fig:NumResSCF} (a) shows the relative error of the
total local potential for each inner SCF iteration, where $x$-axis is
the total number of inner SCF iterations and each lines represent the
inner SCF iterations for an outer SCF iteration. The jump between the
end of previous line and beginning of the next line is introduced by
the update of the nonlocal potential. This is a typical behavior in the
two-level SCF iteration for solving Hartree-Fock-like equations. As the
converged spectral projector in the inner SCF iteration getting closer
to the final convergence, the magnitude of the jump also decreases.
Fig.~\ref{fig:NumResSCF} (b) shows the relative error of the energy
associated with the nonlocal potential for each outer SCF iteration.

\begin{table}[h]
    \centering
    \begin{tabular}{ccccccc}
\toprule
Outer & \multicolumn{3}{c}{GC-ALB} & \multicolumn{3}{c}{Planewave} \\
 SCF & No. $\text{SCF}_{in}$ & $E_{X}$ & rel err
 & No. $\text{SCF}_{in}$ &  $E_{X}$ & rel err\\
\toprule
1 & 10 & -2.825403 & 1.77e-02 & 8 & -2.825400 & 1.77e-02 \\
2 & 7 & -2.841545 & 5.68e-03 & 7 & -2.841543 & 5.68e-03 \\
3 & 6 & -2.849557 & 2.81e-03 & 6 & -2.849555 & 2.81e-03 \\
4 & 5 & -2.852925 & 1.18e-03 & 5 & -2.852922 & 1.18e-03 \\
5 & 6 & -2.854591 & 5.84e-04 & 6 & -2.854588 & 5.84e-04 \\
6 & 5 & -2.855331 & 2.59e-04 & 4 & -2.855328 & 2.59e-04 \\
7 & 6 & -2.855691 & 1.26e-04 & 6 & -2.855688 & 1.26e-04 \\
8 & 5 & -2.855855 & 5.73e-05 & 5 & -2.855852 & 5.73e-05 \\
9 & 5 & -2.855934 & 2.77e-05 & 4 & -2.855932 & 2.80e-05 \\
10 & 4 & -2.855970 & 1.27e-05 & 5 & -2.855968 & 1.24e-05 \\
11 & 3 & -2.855989 & 6.45e-06 & 2 & -2.855986 & 6.49e-06 \\
\bottomrule
    \end{tabular}
    \caption{Comparison of the GC-ALB method and the planewave
    method in self consistent field iteration. The Hamiltonian
    operator defined in \eqref{eqn:HrHF} is solved by a two levels
    of SCF iteration combined with either the GC-ALB method or the
    planewave method. No. $\text{SCF}_{in}$ denotes the number of
    inner SCF iterations, $E_{X}$ denotes the energy associated with
    the nonlocal potential, and rel err is the relative change of
    the $E_{X}$ every outer SCF iteration.} \label{tab:NumResSCF}
\end{table}

\begin{figure}[h]
  \begin{center}
    \subfloat[]{\includegraphics[width=0.45\textwidth]{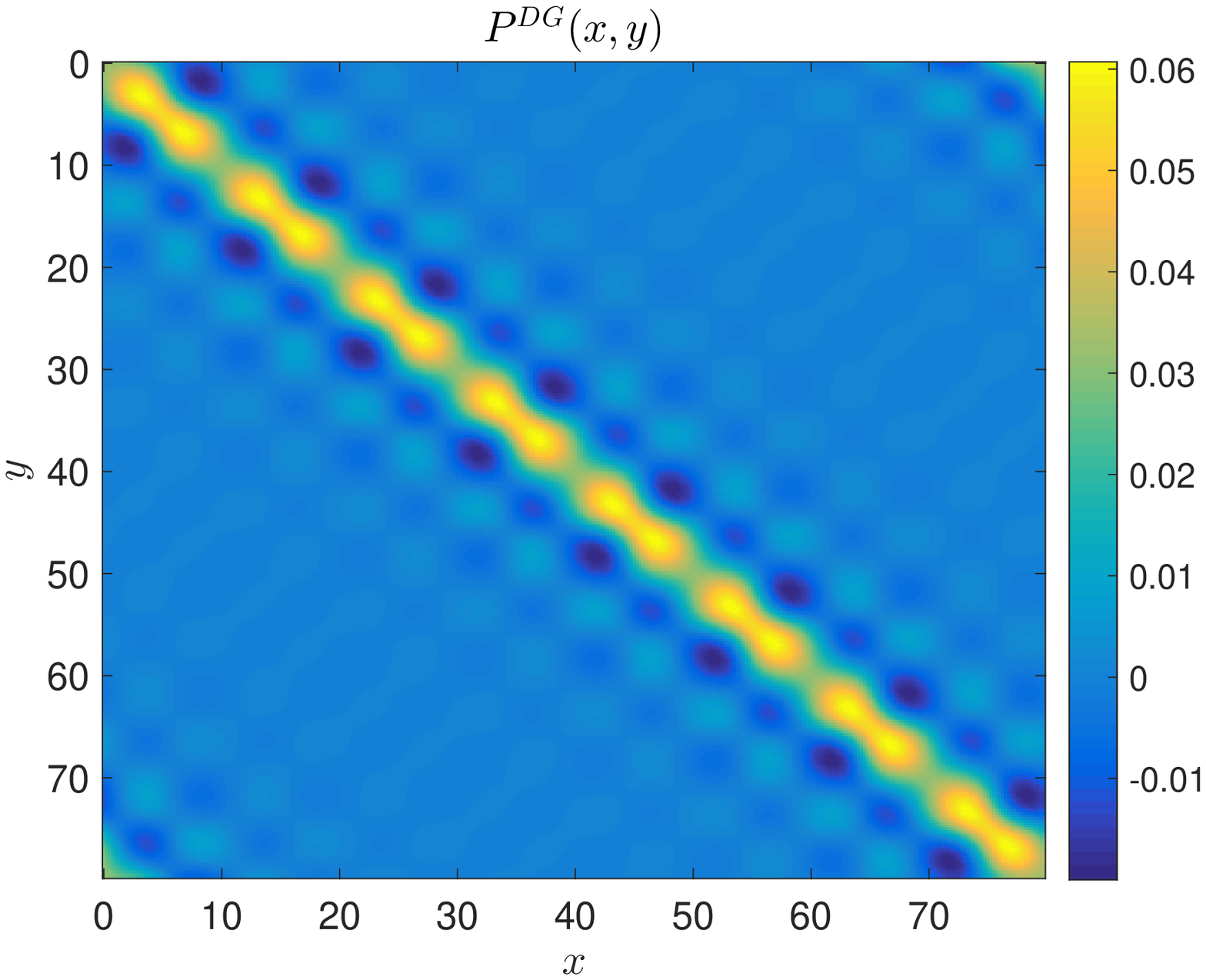}}~~
    \subfloat[]{\includegraphics[width=0.45\textwidth]{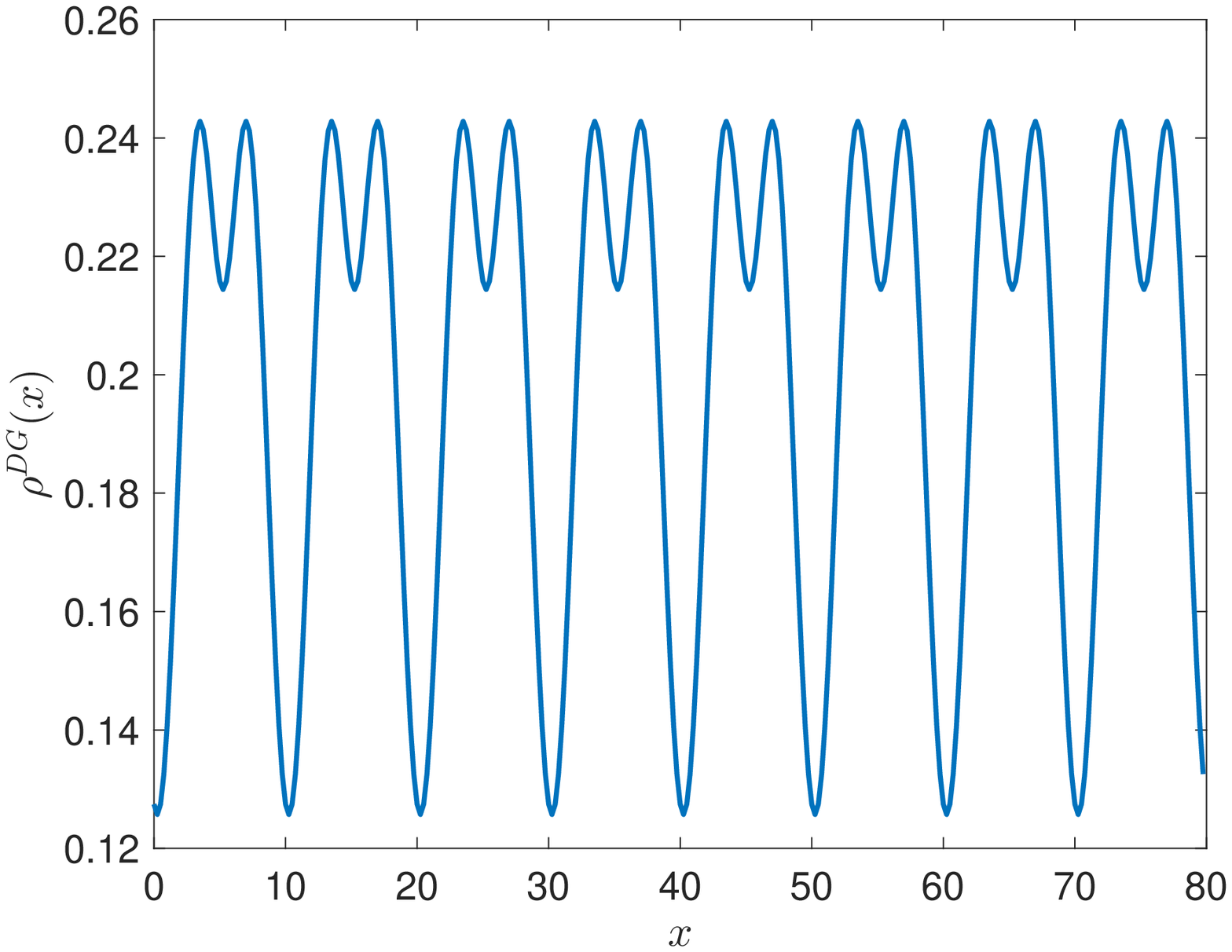}}\\
    \subfloat[]{\includegraphics[width=0.45\textwidth]{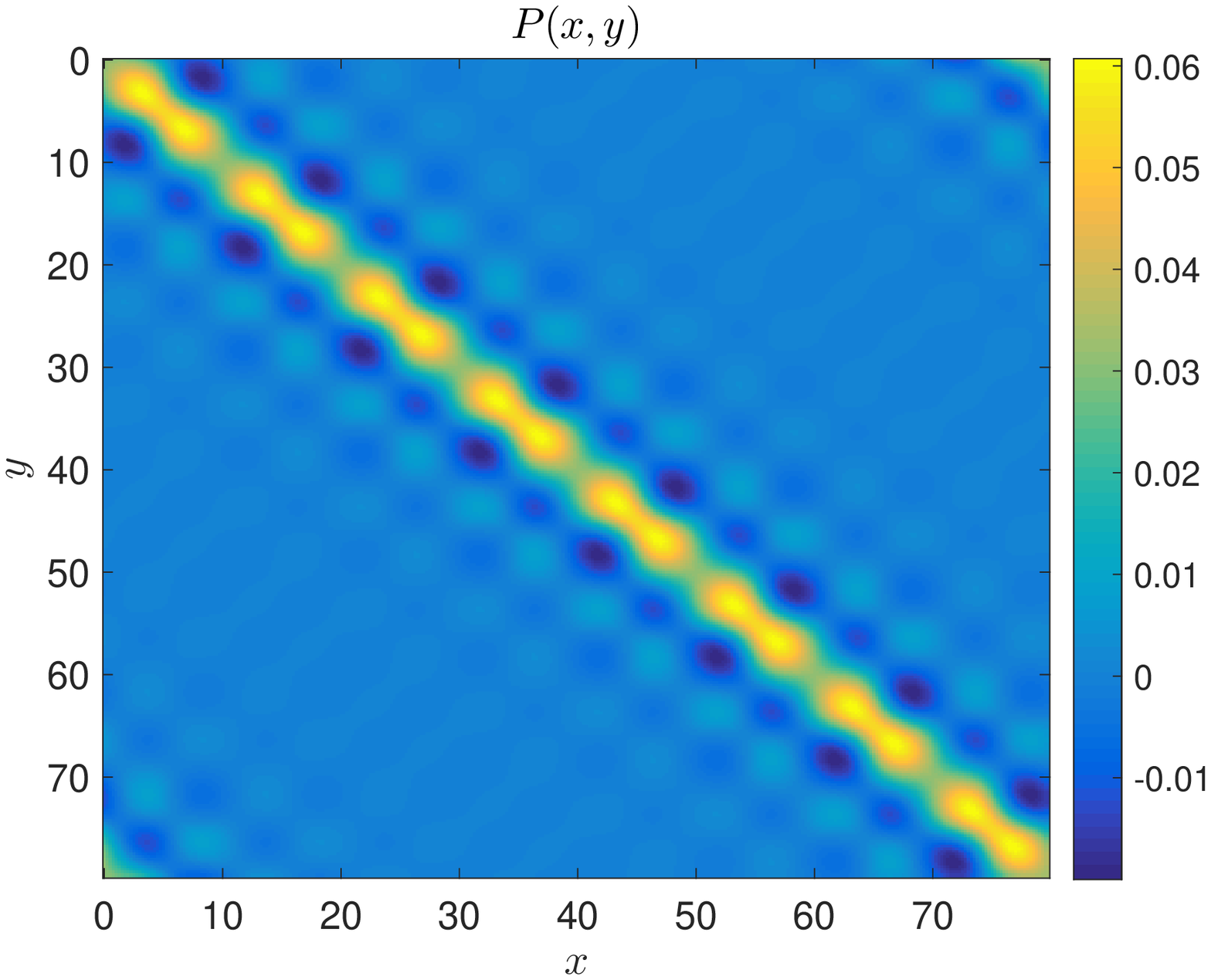}}~~
    \subfloat[]{\includegraphics[width=0.45\textwidth]{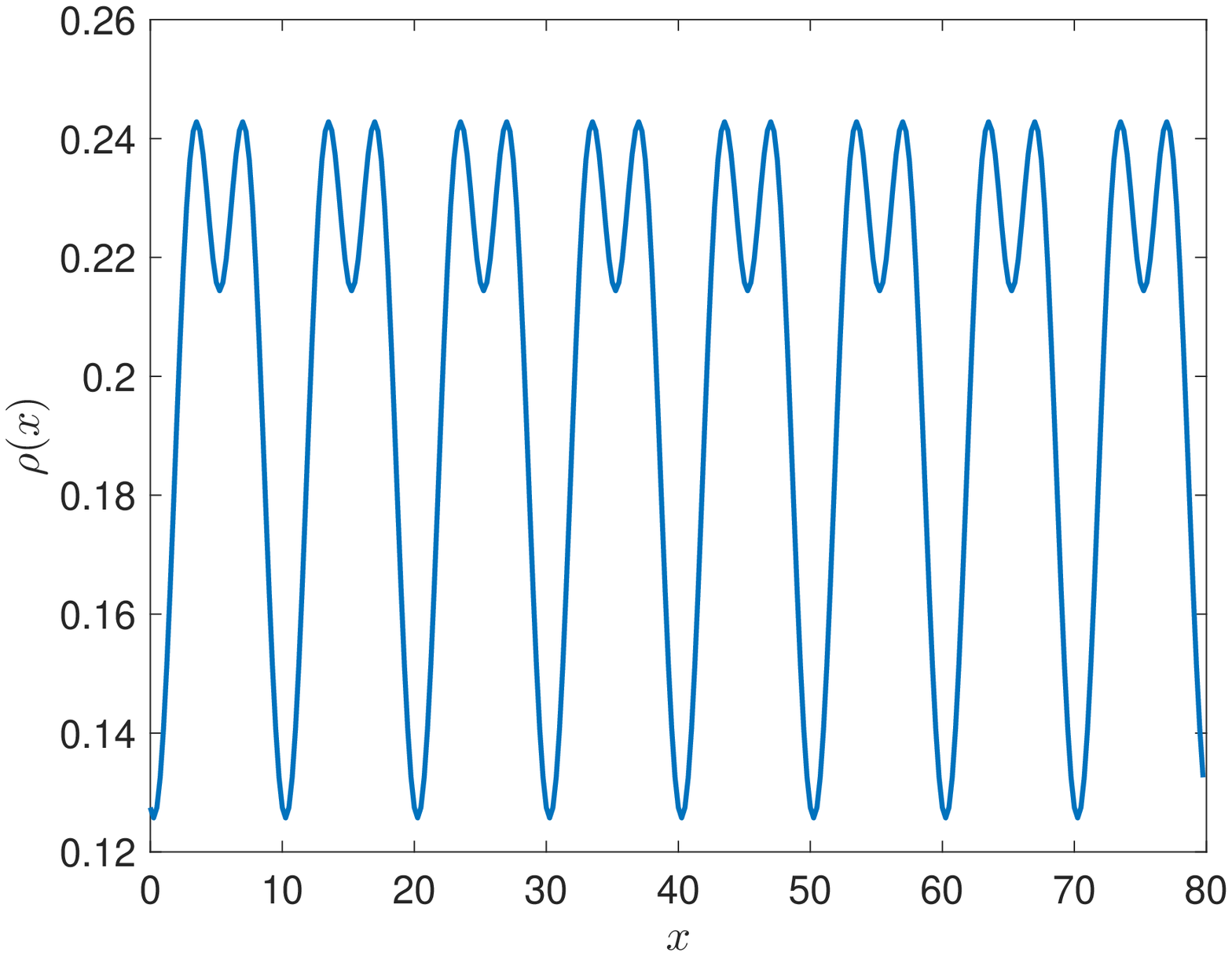}}\\
    \subfloat[]{\includegraphics[width=0.45\textwidth]{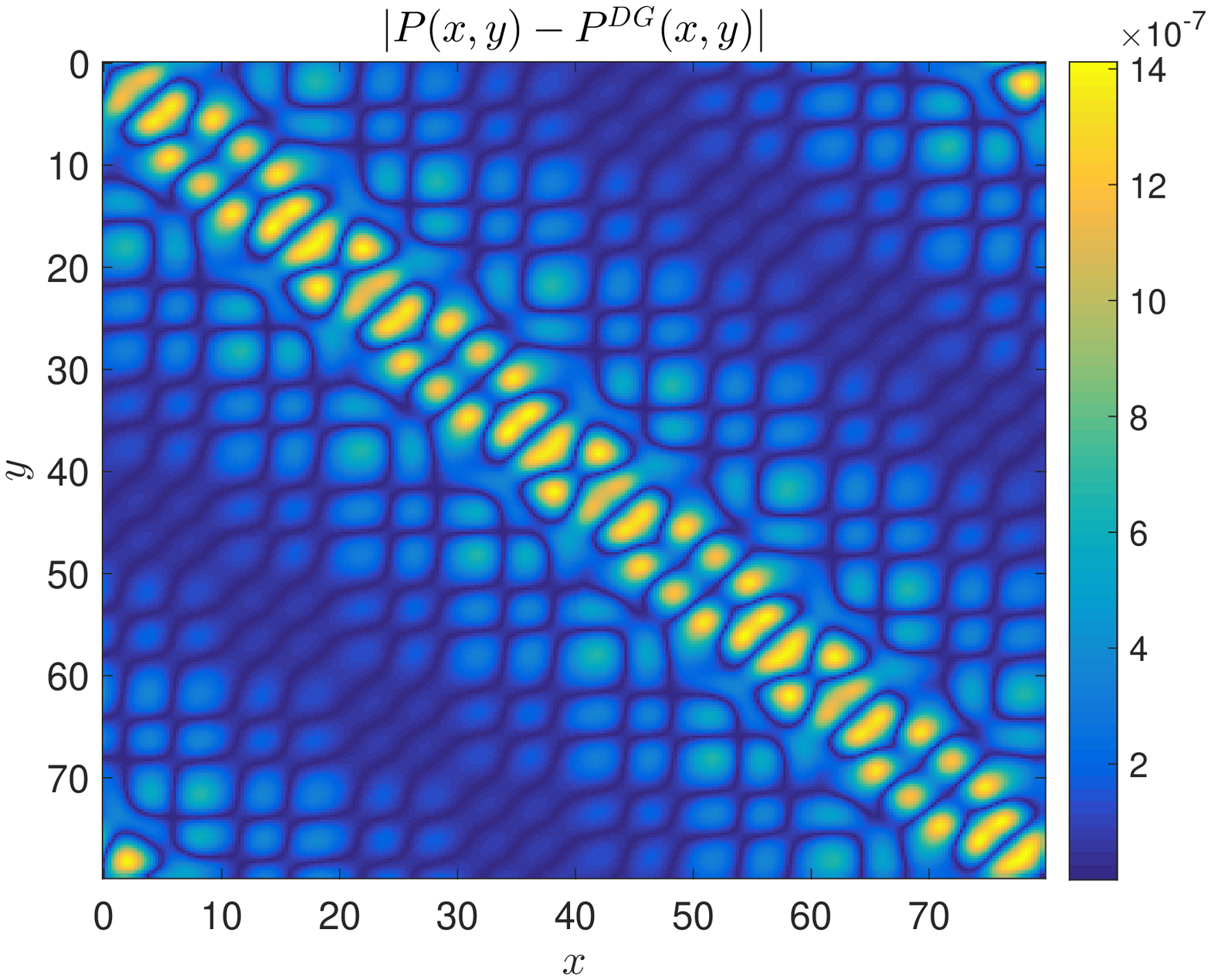}}~~
    \subfloat[]{\includegraphics[width=0.45\textwidth]{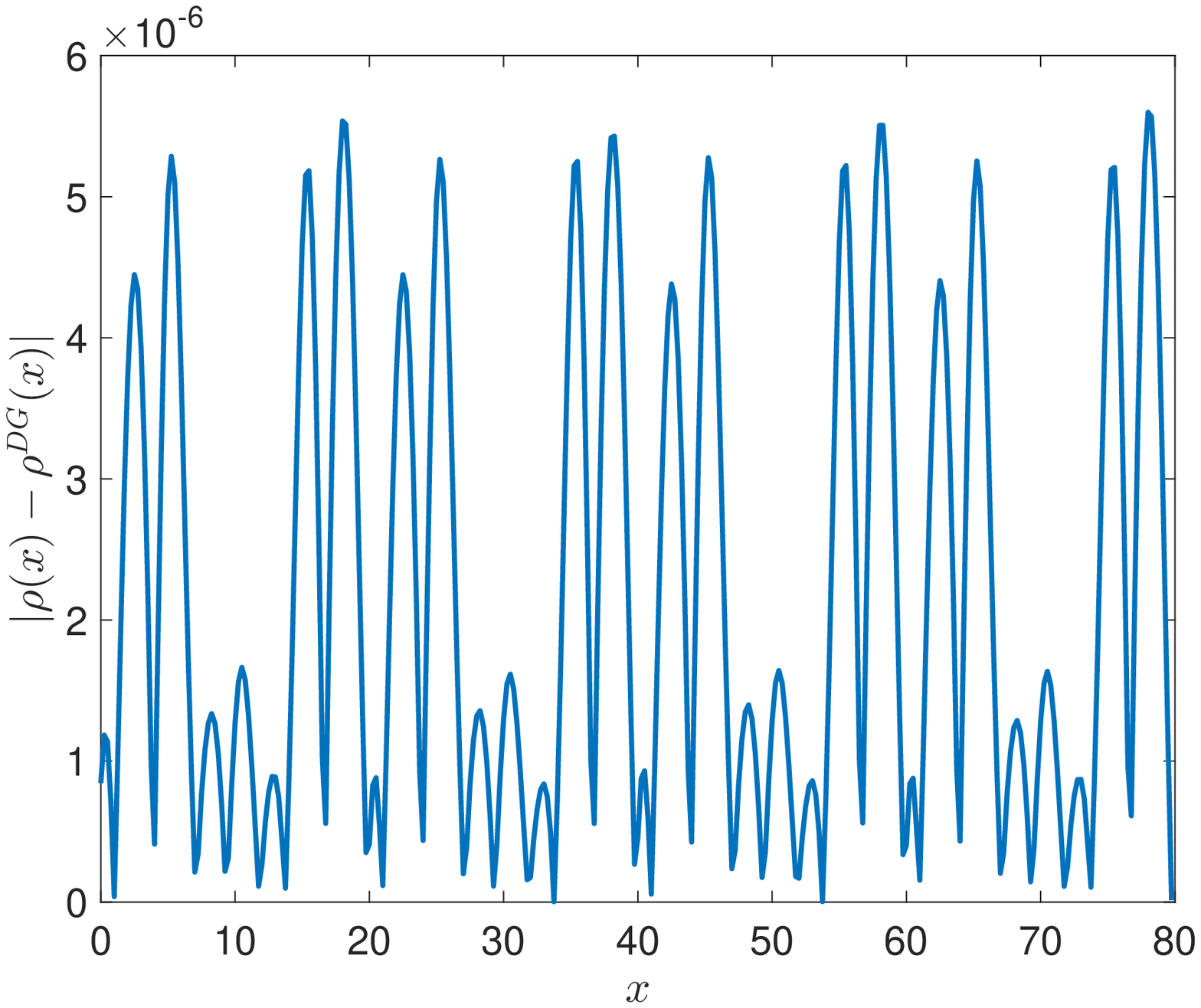}}\\
  \end{center}
  \caption{(a) Kernel of the spectral projector and (b) electron
  density associated with the 1D model for Hartree-Fock-like equation
  calculated by the GC-ALB method, whereas (c) kernel of the spectral
  projector and (d) electron density are calculated by the planewave
  method. (e) is the absolute difference between (a) and (c), and (f)
  is the absolute difference between (b) and (d). } \label{fig:rhoHF}
\end{figure}

Tab.~\ref{tab:NumResSCF} indicates that the calculation using both
the GC-ALB set and the planewave basis set converges within 11 outer
SCF iterations to a relative error around $6\times 10^{-6}$, and the
number of inner iterations in each outer iteration is comparable in
both methods. This indicates that the use of the GC-ALB set does not
increase the number of the SCF iterations in the nonlinear setup. The
relative error from both methods also behaves similarly throughout the
SCF iteration.  The spectral projector, as well as electron density
defined to be diagonal of the converged projector $\rho(x)=P(x,x)$
for both methods are given in Fig.~\ref{fig:rhoHF}. The point-wise
relative differences for the projector and the density are provided at
the last row of Fig.~\ref{fig:rhoHF}, where the errors are about the
same level as that of the relative error in Tab.~\ref{tab:NumResSCF}.

\FloatBarrier

\section{Conclusion}\label{sec:conclusion}

We developed a new method to construct an efficient basis set to
represent the spectral projector of a second order differential
operator $H$ with reduced degrees of freedom.  For a given partition
of the global domain into sub-domains called elements, the optimal
discontinuous basis set on any element can be given by the singular
value decomposition of the matrix row block of the spectral projector
associated with the element. Our globally constructed adaptive local
basis set (GC-ALB) can efficiently approximate such an optimal basis
set in practice. The GC-ALB set can be obtained by only applying a
matrix function $f(H)$ to a small number of random vectors on the
global domain, without the need of any buffer areas to define a
series of local problems. The GC-ALB set can be used in the context
of the discontinuous Galerkin (DG) framework to approximate the
spectral projector on the global domain. When the potential is
local, the reduced DG matrix is a block sparse matrix. Hence the
evaluation of the matrix representation of the spectral projector
can be evaluated using fast methods based on sparse linear algebra
operations, such as the pole expansion and selected inversion
method (PEXSI)~\cite{LinLuYingEtAl2009,LinChenYangEtAl2013},
and the purification methods~\cite{McWeeny1960,Goedecker1999}.
Our method is also flexible and can be applied to operators with
local and nonlocal potentials.  We verified the effectiveness of
the basis set using one, two and three dimensional linear problems,
as well as one-dimensional nonlocal, as well as nonlinear problems
resembling Hartree-Fock problems.  Numerical results indicate that the
GC-ALB set achieve nearly optimal performance in terms of the number
of degrees of freedom per element, which reduces both the storage and
the computational cost. In the near future, we will explore the usage
of the GC-ALB set for Kohn-Sham density functional theory calculations
for real materials.


\section*{Acknowledgments} L. L. was partially supported by the
National Science Foundation under Grant No. DMS-1652330, and by
the U.S.  Department of Energy under contract number DE-SC0017867,
DE-AC02-05CH11231, and the Scientific Discovery through Advanced
Computing (SciDAC) program.  Y. L. was partially supported by the
supported by the National Science Foundation under award DMS-1328230
and the U.S. Department of Energy's Advanced Scientific Computing
Research program under award DE-FC02-13ER26134/DE-SC0009409. We would
like to thank Lexing Ying for fruitful discussions.


\bibliographystyle{siam}
\bibliography{dlrbasis}

\begin{thebibliography}{10}

\bibitem{Akhiezer1990}
{\sc N.~I. Akhiezer}, {\em Elements of the theory of elliptic functions},
  American Mathematical Soc., 1990.

\bibitem{AlemanyJainKronikChelikowsky2004}
{\sc M.~Alemany, M.~Jain, L.~Kronik, and J.~Chelikowsky}, {\em Real-space
  pseudopotential method for computing the electronic properties of periodic
  systems}, Phys. Rev. B, 69 (2004), p.~075101.

\bibitem{Anderson1965}
{\sc D.~G. Anderson}, {\em Iterative procedures for nonlinear integral
  equations}, J. Assoc. Comput. Mach., 12 (1965), pp.~547--560.

\bibitem{Arnold1982}
{\sc D.~N. Arnold}, {\em An interior penalty finite element method with
  discontinuous elements}, SIAM J. Numer. Anal., 19 (1982), pp.~742 -- 760.

\bibitem{ArnoldBrezziCockburnEtAl2002}
{\sc D.~N. Arnold, F.~Brezzi, B.~Cockburn, and L.~D. Marini}, {\em Unified
  analysis of discontinuous {G}alerkin methods for elliptic problems}, SIAM J.
  Numer. Anal., 39 (2002), pp.~1749--1779.

\bibitem{BabuvskaMelenk1997}
{\sc I.~Babu{\v{s}}ka and J.~M. Melenk}, {\em The partition of unity method},
  Int. J. Numer. Meth. Eng., 40 (1997), pp.~727--758.

\bibitem{BanerjeeLinHuEtAl2016}
{\sc A.~S. Banerjee, L.~Lin, W.~Hu, C.~Yang, and J.~E. Pask}, {\em {Chebyshev
  polynomial filtered subspace iteration in the Discontinuous Galerkin method
  for large-scale electronic structure calculations}}, J. Chem. Phys., 145
  (2016), p.~154101.

\bibitem{Becke1993}
{\sc A.~D. Becke}, {\em Density functional thermochemistry. {III}. {T}he role
  of exact exchange}, J. Chem. Phys., 98 (1993), p.~5648.

\bibitem{BerrutTrefethen2004}
{\sc J.-P. Berrut and L.~N. Trefethen}, {\em Barycentric {L}agrange
  interpolation}, SIAM Rev., 46 (2004), pp.~501--517.

\bibitem{BlumGehrkeHankeEtAl2009}
{\sc V.~Blum, R.~Gehrke, F.~Hanke, P.~Havu, V.~Havu, X.~Ren, K.~Reuter, and
  M.~Scheffler}, {\em {Ab initio molecular simulations with numeric
  atom-centered orbitals}}, Comput. Phys. Commun., 180 (2009), pp.~2175--2196.

\bibitem{FattebertBernholc2000}
{\sc J.~L. Fattebert and J.~Bernholc}, {\em {Towards grid-based O(N)
  density-functional theory methods: Optimized nonorthogonal orbitals and
  multigrid acceleration}}, Phys. Rev. B, 62 (2000), pp.~1713--1722.

\bibitem{Garc'ia-CerveraLuXuanEtAl2009}
{\sc C.~J. Garc\'{\i}a-Cervera, Jianfeng Lu, Yulin Xuan, and Weinan E}, {\em
  Linear-scaling subspace-iteration algorithm with optimally localized
  nonorthogonal wave functions for {K}ohn-{S}ham density functional theory},
  Phys. Rev. B, 79 (2009), p.~115110.

\bibitem{GenoveseNeelovGoedeckerEtAl2008}
{\sc L.~Genovese, A.~Neelov, S.~Goedecker, T.~Deutsch, S.~A. Ghasemi,
  A.~Willand, D.~Caliste, O.~Zilberberg, M.~Rayson, A.~Bergman, and
  R.~Schneider}, {\em Daubechies wavelets as a basis set for density functional
  pseudopotential calculations}, J. Chem. Phys., 129 (2008), p.~014109.

\bibitem{GiannozziBaroniBoniniEtAl2009}
{\sc Paolo Giannozzi, Stefano Baroni, Nicola Bonini, Matteo Calandra, Roberto
  Car, Carlo Cavazzoni, Davide Ceresoli, Guido~L Chiarotti, Matteo Cococcioni,
  Ismaila Dabo, Andrea~Dal Corso, Stefano de~Gironcoli, Stefano Fabris, Guido
  Fratesi, Ralph Gebauer, Uwe Gerstmann, Christos Gougoussis, Anton Kokalj,
  Michele Lazzeri, Layla Martin-Samos, Nicola Marzari, Francesco Mauri,
  Riccardo Mazzarello, Stefano Paolini, Alfredo Pasquarello, Lorenzo Paulatto,
  Carlo Sbraccia, Sandro Scandolo, Gabriele Sclauzero, Ari~P Seitsonen,
  Alexander Smogunov, Paolo Umari, and Renata~M Wentzcovitch}, {\em {QUANTUM
  ESPRESSO}: {A} modular and open-source software project for quantum
  simulations of materials}, J. Phys.: Condens. Matter, 21 (2009),
  pp.~395502--395520.

\bibitem{Goedecker1999}
{\sc S.~Goedecker}, {\em {Linear scaling electronic structure methods}}, Rev.
  Mod. Phys., 71 (1999), pp.~1085--1123.

\bibitem{GolubVan2013}
{\sc G.~H. Golub and C.~F. Van~Loan}, {\em Matrix computations}, Johns Hopkins
  Univ. Press, Baltimore, fourth~ed., 2013.

\bibitem{Guttel2015a}
{\sc S.~G{\"{u}}ttel, E.~Polizzi, P.~T.~P. Tang, and G.~Viaud}, {\em
  {{Zolotarev} quadrature rules and load balancing for the {FEAST}
  eigensolver}}, SIAM J. Sci. Comput., 37 (2015), pp.~A2100--A2122.

\bibitem{HalkoMartinssonTropp2011}
{\sc N.~Halko, P.-G. Martinsson, and J.~A. Tropp}, {\em Finding structure with
  randomness: Probabilistic algorithms for constructing approximate matrix
  decompositions}, SIAM Rev., 53 (2011), pp.~217--288.

\bibitem{HeydScuseriaErnzerhof2003}
{\sc J.~Heyd, G.~E. Scuseria, and M.~Ernzerhof}, {\em Hybrid functionals based
  on a screened coulomb potential}, J. Chem. Phys., 118 (2003), pp.~8207--8215.

\bibitem{HohenbergKohn1964}
{\sc P.~Hohenberg and W.~Kohn}, {\em {Inhomogeneous electron gas}}, Phys. Rev.,
  136 (1964), pp.~B864--B871.

\bibitem{HuLinYang2015a}
{\sc W.~Hu, L.~Lin, and C.~Yang}, {\em {DGDFT}: A massively parallel method for
  large scale density functional theory calculations}, J. Chem. Phys., 143
  (2015), p.~124110.

\bibitem{Knyazev2001}
{\sc A.~V. Knyazev}, {\em Toward the optimal preconditioned eigensolver:
  Locally optimal block preconditioned conjugate gradient method}, SIAM J. Sci.
  Comp., 23 (2001), pp.~517--541.

\bibitem{KohnSham1965}
{\sc W.~Kohn and L.~Sham}, {\em {Self-consistent equations including exchange
  and correlation effects}}, Phys. Rev., 140 (1965), pp.~A1133--A1138.

\bibitem{Lanczos1950}
{\sc C.~Lanczos}, {\em An iteration method for the solution of the eigenvalue
  problem of linear differential and integral operators}, J. Res. Nat. Bur.
  Stand., 45 (1950), pp.~255--282.

\bibitem{Li2017b}
{\sc Y.~Li and H.~Yang}, {\em {Spectrum slicing for sparse {Hermitian} definite
  matrices based on {Zolotarev's} functions}}, tech. report, 2017.

\bibitem{Lin2016ACE}
{\sc L.~Lin}, {\em Adaptively compressed exchange operator}, J. Chem. Theory
  Comput., 12 (2016), p.~2242.

\bibitem{LinChenYangEtAl2013}
{\sc L.~Lin, M.~Chen, C.~Yang, and L.~He}, {\em Accelerating atomic
  orbital-based electronic structure calculation via pole expansion and
  selected inversion}, J. Phys.: Condens. Matter, 25 (2013), p.~295501.

\bibitem{LinLuYingEtAl2009}
{\sc L.~Lin, J.~Lu, L.~Ying, R.~Car, and W.~E}, {\em Fast algorithm for
  extracting the diagonal of the inverse matrix with application to the
  electronic structure analysis of metallic systems}, Comm. Math. Sci., 7
  (2009), p.~755.

\bibitem{LinLuYingE2012}
{\sc L.~Lin, J.~Lu, L.~Ying, and W.~E}, {\em {Adaptive local basis set for
  Kohn-Sham density functional theory in a discontinuous Galerkin framework I:
  Total energy calculation}}, J. Comput. Phys., 231 (2012), pp.~2140--2154.

\bibitem{LinLuYingE2012a}
\leavevmode\vrule height 2pt depth -1.6pt width 23pt, {\em Optimized local
  basis function for {Kohn-Sham} density functional theory}, J. Chem. Phys.,
  231 (2012), p.~4515.

\bibitem{LinStamm2016}
{\sc L.~Lin and B.~Stamm}, {\em A posteriori error estimates for discontinuous
  {G}alerkin methods using non-polynomial basis functions. {P}art {I}: {S}econd
  order linear {PDE}}, Math. Model. Numer. Anal., 50 (2016), p.~1193.

\bibitem{LinStamm2017}
\leavevmode\vrule height 2pt depth -1.6pt width 23pt, {\em A posteriori error
  estimates for discontinuous {G}alerkin methods using non-polynomial basis
  functions. {P}art {II}: {E}igenvalue problems}, Math. Model. Numer. Anal.,
  (2017, in press).

\bibitem{LinYang2013}
{\sc L.~Lin and C.~Yang}, {\em {Elliptic preconditioner for accelerating self
  consistent field iteration in Kohn-Sham density functional theory}}, SIAM J.
  Sci. Comp., 35 (2013), pp.~S277--S298.

\bibitem{Martin2004}
{\sc R.~Martin}, {\em Electronic structure -- Basic theory and practical
  methods}, Cambridge Univ. Pr., West Nyack, {NY}, 2004.

\bibitem{McWeeny1960}
{\sc R.~McWeeny}, {\em Some recent advances in density matrix theory}, Rev.
  Mod. Phys., 32 (1960), pp.~335--369.

\bibitem{MohrRatcliffBoulangerEtAl2014}
{\sc S.~Mohr, L.~E. Ratcliff, P.~Boulanger, L.~Genovese, D.~Caliste,
  T.~Deutsch, and S.~Goedecker}, {\em Daubechies wavelets for linear scaling
  density functional theory}, J. Chem. Phys., 140 (2014), p.~204110.

\bibitem{MotamarriGavini2014}
{\sc P.~Motamarri and V.~Gavini}, {\em Subquadratic-scaling subspace projection
  method for large-scale {Kohn-Sham} density functional theory calculations
  using spectral finite-element discretization}, Phys. Rev. B, 90 (2014),
  p.~115127.

\bibitem{Ozaki:03}
{\sc T.~Ozaki}, {\em Variationally optimized atomic orbitals for large-scale
  electronic structures}, Phys. Rev. B, 67 (2003), pp.~155108--155112.

\bibitem{PaigeSaunders1975}
{\sc C.~C. Paige and M.~A Saunders}, {\em Solution of sparse indefinite systems
  of linear equations}, SIAM J. Numer. Anal., 12 (1975), pp.~617--629.

\bibitem{PaskKleinFongEtAl1999}
{\sc J.~E. Pask, B.~M. Klein, C.~Y. Fong, and P.~A. Sterne}, {\em Real-space
  local polynomial basis for solid-state electronic-structure calculations: {A}
  finite-element approach}, Phys. Rev. B, 59 (1999), p.~12352.

\bibitem{RaysonBriddon2009}
{\sc M.~J. Rayson and P.~R. Briddon}, {\em Highly efficient method for
  {Kohn-Sham} density functional calculations of 500--10000 atom systems},
  Phys. Rev. B, 80 (2009), p.~205104.

\bibitem{SaadSchultz1986}
{\sc Y.~Saad and M.~H. Schultz}, {\em {GMRES}: {A} generalized minimal residual
  algorithm for solving nonsymmetric linear systems}, SIAM J. Sci. Stat.
  Comput., 7 (1986), pp.~856--869.

\bibitem{SchwarzBlahaMadsen2002}
{\sc K.~Schwarz, P.~Blaha, and G.~K.~H. Madsen}, {\em {Electronic structure
  calculations of solids using the WIEN2k package for material sciences}},
  Comput. Phys. Commun., 147 (2002), pp.~71--76.

\bibitem{SkylarisHaynesMostofiEtAl2005}
{\sc C.K. Skylaris, P.D. Haynes, A.A. Mostofi, and M.C. Payne}, {\em
  {Introducing ONETEP: Linear-scaling density functional simulations on
  parallel computers}}, J. Chem. Phys., 122 (2005), p.~084119.

\bibitem{SolerArtachoGaleEtAl2002}
{\sc J.~M. Soler, E.~Artacho, J.~D. Gale, A.~Garc{\'\i}a, J.~Junquera,
  P.~Ordej{\'o}n, and D.~S{\'a}nchez-Portal}, {\em {The SIESTA method for ab
  initio order-N materials simulation}}, J. Phys.: Condens. Matter, 14 (2002),
  pp.~2745--2779.

\bibitem{SukumarPask2009}
{\sc N.~Sukumar and J.~E. Pask}, {\em Classical and enriched finite element
  formulations for {B}loch-periodic boundary conditions}, Int. J. Numer. Meth.
  Engng., 77 (2009), p.~1121.

\bibitem{SzaboOstlund1989}
{\sc A.~Szabo and N.S. Ostlund}, {\em Modern quantum chemistry: Introduction to
  advanced electronic structure theory}, McGraw-Hill, New York, 1989.

\bibitem{GijzenErlanggaVuik2007}
{\sc M.~B. van Gijzen, Y.~A. Erlangga, and C.~Vuik}, {\em Spectral analysis of
  the discrete {H}elmholtz operator preconditioned with a shifted {L}aplacian},
  SIAM J. Sci. Comput., 29 (2007), pp.~1942--1958.

\bibitem{ZhangLinHuEtAl2017}
{\sc G.~Zhang, L.~Lin, W.~Hu, C.~Yang, and J.~E. Pask}, {\em {Adaptive local
  basis set for Kohn--Sham density functional theory in a discontinuous
  Galerkin framework II: Force, vibration, and molecular dynamics
  calculations}}, J. Comput. Phys., 335 (2017), p.~426.

\bibitem{ZhouChelikowskySaad2014}
{\sc Y.~Zhou, J.~R. Chelikowsky, and Y.~Saad}, {\em {Chebyshev-filtered
  subspace iteration method free of sparse diagonalization for solving the
  Kohn-Sham equation}}, J. Comput. Phys., 274 (2014), pp.~770--782.

\bibitem{ZhouSaadTiagoEtAl2006}
{\sc Y.~Zhou, Y.~Saad, M.~L. Tiago, and J.~R. Chelikowsky}, {\em
  Self-consistent-field calculations using {C}hebyshev-filtered subspace
  iteration}, J. Comput. Phys., 219 (2006), pp.~172--184.

\bibitem{Zolotarev1877}
{\sc E.~Zolotarev}, {\em {Application of elliptic functions to questions of
  functions deviating least and most from zero}}, Zap. Imp. Akad. Nauk. St.
  Petersbg., 30 (1877), pp.~1--59.

\end{thebibliography}

\end{document}